\documentclass[ejs,preprint]{imsart}
\usepackage[utf8]{inputenc}
\usepackage[T1]{fontenc}
\usepackage[english]{babel}

\usepackage{amsmath,amsfonts,amssymb,amsthm,bbm,amstext,amsbsy}
\usepackage{graphicx}
\usepackage{natbib}
\usepackage{color}

\usepackage[center,sc,small]{subfigure}

\usepackage{verbatim}

\startlocaldefs

\newcommand{\1}{{\mathbbm 1}}
\newcommand{\E}{{\mathbb E}}
\renewcommand{\P}{{\mathbb P}}
\newcommand{\mR}{{\mathbb R}}
\newcommand{\G}{{\mathbf G}}

\newtheorem{Rem}{Remark}

\newtheorem{Lem}{Lemma}[section]

\newenvironment{Proof} {\noindent {\textbf{Proof}}} { \hfill $\Box$ \\ }
\newtheorem{Prop}{Proposition}[section]

\endlocaldefs

\begin{document}

\begin{frontmatter}
\title{Uniform convergence and asymptotic confidence bands for model-assisted estimators of the mean of sampled functional data}
\runtitle{Model-assisted estimators for functional data}

\begin{aug}
  \author{\fnms{Herv\'e}  \snm{Cardot}\ead[label=hc]{herve.cardot(AT)u-bourgogne.fr}},
  \and
  \author{\fnms{Camelia} \snm{Goga}\ead[label=cg]{camelia.goga(AT)u-bourgogne.fr}}
    \address{
UMR 5584 CNRS Institut de Math\'ematiques de Bourgogne,\\
Universit\'e de Bourgogne, UFR Sciences et Techniques,\\
9 avenue Alain Savary -- BP 47870, 
21078 Dijon Cedex, France.\\
\printead{hc,cg,pl}}

  \author{\fnms{Pauline}  \snm{Lardin}\corref{}%
  \ead[label=pl]{pauline.lardin(AT)laposte.fr}}%

  \runauthor{H.~Cardot, C.~Goga and P. Lardin}

  \affiliation{Universit\'e de Bourgogne}

  \address{
EDF, R{\&}D, ICAME-SOAD, 1 av. du G\'en\'eral de Gaulle, \\ 92141 Clamart, FRANCE \\
\printead{pl}} 


\end{aug}


\begin{abstract}

When the study variable is functional and storage capacities are limited or transmission costs are high, selecting with survey sampling techniques a small fraction of the observations is an interesting alternative to signal compression techniques, particularly when the goal is the estimation of  simple quantities such  as means or totals. We extend, in this functional framework,  model-assisted estimators with linear regression models that can take account of auxiliary variables whose  totals over the population are known. 
We first show, under weak hypotheses on the sampling design and the regularity of the trajectories,  that the estimator  of the mean function as well as its variance estimator are uniformly consistent. Then, under additional assumptions, we prove a functional central limit theorem and we assess rigorously   a fast technique based on simulations of Gaussian processes which is employed  to build asymptotic confidence bands.  The accuracy of the variance function estimator is evaluated on a real dataset of sampled electricity consumption curves measured every half an hour over a period of one week.
\end{abstract}



\end{frontmatter}

\section{Introduction}

\begin{figure}
\begin{center}
\includegraphics[width=\textwidth]{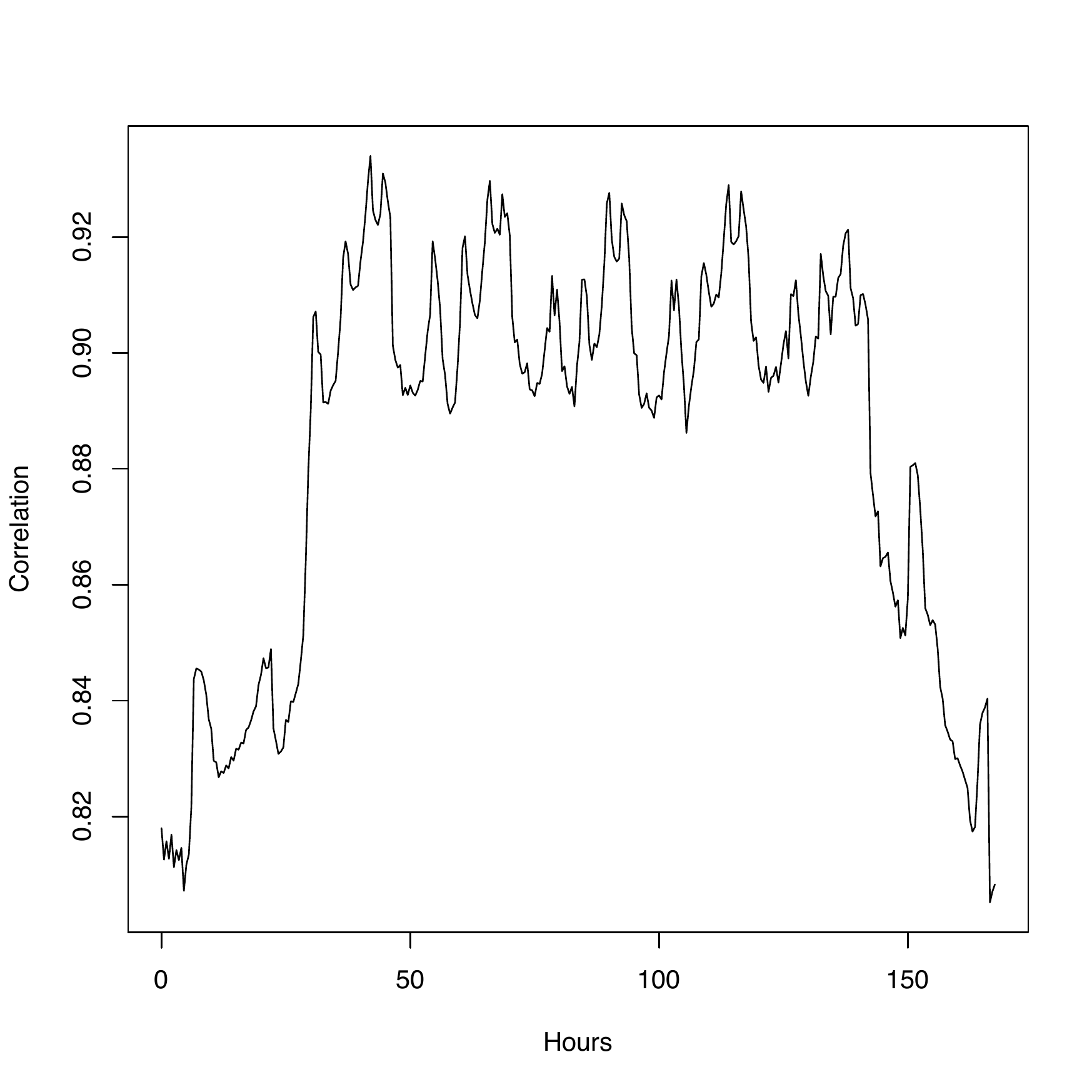}
\caption{Correlation between  the current consumption at each instant $t$ of the week under study and the total past consumption of the week before.}  
\label{fig:tempcorr}
\end{center}
\end{figure}

Survey sampling techniques, which consist in randomly selecting only a part of the elements of a population, are interesting alternatives to signal compression when one has to deal with very large populations of  quantities that evolve  along time. With the development of automatic sensors such very large datasets of temporal data are not unusual anymore  and survey sampling techniques  offer a good trade-off between accuracy of the estimators and size of the analyzed data. Examples  can be found in different domains such as internet traffic monitoring (see \cite{callado2009}) or estimation of energy consumption measured by individual  smart meters. Motivated by the estimation of mean consumption electricity profiles measured every half an hour over one week, \citet{cardot_josserand_2011} have introduced  Horvitz-Thompson estimators of the mean function and have shown, under weak hypotheses on the regularity of the functional trajectories and the sampling design, that one gets  uniformly convergent estimators.  They also prove a functional central limit theorem, in the space of continuous functions, that can, in part, justify the construction of asymptotic confidence bands.  
More recently, \cite{CDGJL2012} made a  comparison,  in terms of precision of the mean estimators of electricity load curves and width of the confidence bands, of different sampling approaches that can take auxiliary information into account.
One of the conclusions of this empirical study was  that  very simple strategies based on   simple sampling designs (such as simple random sampling without replacement) could be improved much if some well chosen auxiliary information, whose total is known for the whole population,   is also taken into account at the estimation stage, with model-assisted estimators.
 Important variables for the electricity consumption such as temperature or geographical location were not available for these datasets so that only one  auxiliary information, the mean past consumption over the previous period, was taken into account. Its correlation with the current consumption is always very high (see Figure~\ref{fig:tempcorr}) so that linear regression models are natural candidates for assisting the Horvitz-Thompson estimator. More generally, one advantage of linear approaches is 
that they only require the knowledge of the auxiliary variable  totals in the population. More sophisticated nonlinear or nonparametric approaches would have required to know the values of the auxiliary variables for all the elements of the population.

Thus, we focus in this paper on linear relationships between the set of auxiliary variables and the response at each instant $t$ of the current period. The regression coefficients vary in time (see \citet{faraway_1997} or \cite{Ramsay_Silverman_Livre})  
so that the model-assisted estimator can be seen as a direct extension, to a functional or varying-time context, of the generalized regression  (GREG) estimators studied in \cite{RobinsonSarndal83} and \cite{SarndalLivre}. Note also that from another point of view,  the model-assisted estimator can be obtained using a calibration technique (\cite{DevilleSarndal1992}).

Confidence bands are then built using a simulation technique  developed in  \cite{faraway_1997},  \cite{CuevasFF06} and \cite{degras_2011}. We first estimate  the covariance function of the mean estimator and then,  assuming asymptotic normality, perform simulations of a centered Gaussian process whose covariance function is the covariance function estimated at the previous step. We can, this way, obtain an approximation to the law of the "sup" and deduce confidence bands for the mean trajectory. In a recent work, \cite{CardotDegrasJosserand2012} have given a rigorous mathematical justification of this technique for sampled functional data and Horvitz-Thompson estimators for the mean. The required theoretical ingredients that can justify  such a procedure are the functional central limit theorem for the mean estimator, in the space of continuous functions equipped with the sup-norm, as well as a  uniformly consistent estimator of the variance function. 

The aim of this paper is to study the asymptotic properties of model-assisted estimators and to show that we obtain, under classical assumptions, a uniformly consistent estimator of the mean as well as of its variance function. One additional difficulty is that, for model-assisted estimators, the variance function cannot be derived exactly and we can only have asymptotic approximations. Then, we  deduce that the confidence bands built via simulations have asymptotically the desired coverage. 
In Section 2, we introduce notations and we suggest a slight modification of  the model-assisted estimators which permits control of the variance of  the regression coefficient estimator. Under classical assumptions on the sampling design and on the regularity of the trajectories, we state, in Section 3,   the uniform convergence of the model assisted-estimators to the mean function. Under additional assumptions on the design we also prove that we can get a consistent estimator of the covariance function and a functional central limit theorem that can justify rigorously that the confidence bands built with the procedure  based on Gaussian process simulations attain asymptotically the desired level of confidence. In Section~4, we assess the precision of the variance estimator on the real dataset consisting of electricity consumption curves studied in \cite{CDGJL2012} and observe that, in our context, the approximation error is negligible compared to the sampling error. 
A brief discussion about possible extensions and future investigation is proposed in Section~5. All the proofs are gathered in an Appendix.
  
\section{Notations and estimators}

\subsection{The Horvitz Thompson estimator for functional data}
Let us consider a finite population $U_N=\{1,...,N\}$ of size $N$ supposed to be known, and suppose that, for each unit $k$ of the population $U_N$, we can observe a deterministic curve $Y_k=(Y_k(t))_{t\in[0,T]}$.
The target is   the mean trajectory $\mu_N(t)$, $t\in[0,T]$, defined as follows:
\begin{align}
\mu_N(t)&=\frac{1}{N}\sum_{k\in U}Y_k(t).
\end{align}

We consider a sample $s$, with size $n$,  drawn from $U_N$ according to a fixed-size sampling design $p_N(s)$, where $p_N(s)$ is the probability of drawing the sample $s$. For simplicity of notations, the subscript $N$ is omitted when there is no ambiguity. 
 We suppose that the first and second order inclusion probabilities satisfy $\pi_k= \P(k\in s)>0$, for all $k\in U$, and $\pi_{kl}= \P(k\&l\in s)>0$ for all $k,l\in U_N$, $k\neq l$. 
Without auxiliary information, the population mean curve $\mu(t)$ is often estimated by the Horvitz-Thompson estimator, defined as follows for $t \in [0,T],$ 
\begin{align}
\widehat{\mu}(t)&=\frac{1}{N}\sum_{k\in s}\frac{Y_k(t)}{\pi_k}=\frac{1}{N}\sum_{k\in U}\frac{Y_k(t)}{\pi_k} \1_{k }, 
\label{horvitz-thompson}
\end{align}
where $\1_{k}$ is the sample membership indicator, $\1_{k }=1$ if $k\in s$ and $\1_{k}=0$ otherwise. For each $t\in [0,T],$  the estimator $\widehat{\mu}(t)$ is design-unbiased for $\mu(t)$, $i.e.$ $\E_p(\widehat{\mu}(t))={\mu}(t)$, where $\E_p[.]$ denotes expectation with respect to the sampling design.

The Horvitz-Thompson covariance function of $\widehat{\mu}$ between two instants $r$ and $t,$ computed with respect to the sampling design, is defined as follows
\begin{align}
\mbox{Cov}_p(\widehat{\mu}(r),\widehat{\mu}(t)) &=\frac{1}{N^2}\sum_{k \in U}\sum_{l \in U}(\pi_{kl}-\pi_k\pi_l)\frac{Y_k(r)}{\pi_k}\cdot\frac{Y_l(t)}{\pi_l} \quad r,t\in [0,T].\label{covar_ht}
\end{align} 
Note that for $r=t,$ we obtain the Horvitz-Thompson variance function.

\subsection{The  mean curve  estimator assisted by a functional linear model}
Let us  suppose now that for each unit $k\in U$ we can also observe $p$ real variables, $X_1,...,X_p,$  and let us denote by  $\mathbf{x}_k=(x_{k1},...,x_{kp})'$, the value of the auxiliary variable vector for each unit $k$ in the population.
We introduce an estimator based on a linear regression model that can use these variables in order to improve the accuracy of $\widehat{\mu}.$
By analogy to the real case (see \textit{e.g.} \cite{SarndalLivre}) we suppose that the relationship between the functional variable of interest and the auxiliary variables is modeled by the superpopulation model $\xi$ defined as follows:
\begin{align}
\xi: \quad Y_k(t)&= \mathbf{x}_k'\boldsymbol{\beta}(t)+\epsilon_{kt}, \quad t\in [0,T]
\label{model:xi}
\end{align}
where $\boldsymbol{\beta}(t) = (\beta_1(t), \ldots, \beta_p(t))'$ is the vector of functional regression coefficients, $\epsilon_{kt}$ are independent (across units) and  centered continuous time processes, $\E_{\xi}(\epsilon_{kt})=0,$ with covariance function $\mbox{Cov}_{\xi}(\epsilon_{kt},\epsilon_{kr}) = \Gamma(t,r),$ for $(t,r) \in [0,T]\times[0,T].$ 
This model is a direct extension to several variables of the functional linear model proposed by \citet{faraway_1997}.

If $\mathbf{x}_k$ and $Y_k$ are known for all units $k\in U$ and if the matrix $\mathbf{G}=\frac{1}{N}\sum_{k\in U}\mathbf{x}_k\mathbf{x}_k'$ is invertible, it is possible, under the model $\xi,$ to estimate $\boldsymbol\beta(t)$ by $\tilde{\boldsymbol{\beta}}(t)=\mathbf{G}^{-1}\frac{1}{N}\sum_{k \in U}\mathbf{x}_kY_k(t),$ the ordinary least squares estimator. Then, the  mean curve $\mu(t)$ can be estimated by the generalized difference estimator (see \citet{SarndalLivre}, Chapter 6) defined as follows for all $t \in [0,T],$
\begin{align}
\tilde{\mu}(t) & =  \frac{1}{N}\sum_{k\in U}\mathbf{x}_k'\tilde{\boldsymbol\beta}(t)-\frac{1}{N}\sum_{k\in s}\frac{\mathbf{x}_k'\tilde{\boldsymbol\beta}(t)-Y_k(t)}{\pi_k}\label{mu_difference}\\
  & = \frac{1}{N}\sum_{k\in U}\tilde{Y}_k(t)-\frac{1}{N}\sum_{k\in s}\frac{\tilde{Y}_k(t)-Y_k(t)}{\pi_k}, 
  \nonumber
  \end{align}
where $\tilde{Y}_k(t) = \mathbf{x}_k'\tilde{\boldsymbol\beta}(t).$

In practice, we do not know $Y_k$ except for $k\in s,$  and it is not possible to compute $\tilde{\boldsymbol\beta}(t)$. An estimator of $\mu(t)$ is obtained by substituting  each total in $\tilde{\boldsymbol\beta}(t)$ by its Horvitz-Thompson estimator. Thus, if the matrix $\widehat{\mathbf{G}}=\frac{1}{N}\sum_{k\in s}\frac{\mathbf{x}_k\mathbf{x}_k'}{\pi_k}$ is invertible, $\tilde{\boldsymbol\beta}(t)$ is estimated by:
\begin{align*}
\widehat{\boldsymbol\beta}(t) & =\widehat{\mathbf{G}}^{-1}\frac{1}{N}\sum_{k \in s}\frac{\mathbf{x}_kY_k(t)}{\pi_k}, \quad t\in[0,T].
\end{align*}
\noindent Remark that the denominator $N$ is used in the expression of $\tilde{\boldsymbol\beta}(t)$ for asymptotic purposes  and need not be estimated.  The model-assisted estimator $\widehat{\mu}_{MA}(t)$ is then defined by replacing $\tilde{\boldsymbol\beta}(t)$ by $\widehat{\boldsymbol\beta}(t)$ in (\ref{mu_difference}), 
 \begin{align}
\widehat{\mu}_{MA}(t)&= \frac{1}{N}\sum_{k\in U}\widehat{Y}_k(t)-\frac{1}{N}\sum_{k\in s}\frac{\widehat{Y}_k(t)-Y_k(t)}{\pi_k}, \quad t\in[0,T],
\label{mu:ma}
 \end{align}
where $\widehat{Y}_k(t)=\mathbf{x}_k'\widehat{\boldsymbol\beta}(t).$ 
Since $\sum_{k\in U}\widehat{Y}_k(t)= \left( \sum_{k\in U} \mathbf{x}_k \right)'\widehat{\boldsymbol\beta}(t)$, the only required information to build  $\widehat{\mu}_{MA}(t)$ is $\mathbf{x}_k$ and $Y_k(t)$ for all the units $k \in s$ as well as the population totals of the auxiliary variables, $\sum_{k \in U} \mathbf{x}_k$.

\begin{Rem}
If the vector of auxiliary information contains the intercept (constant term), then it can be  shown (see  \cite{Sarndal1980}) that the Horvitz-Thompson estimator of the estimated residuals $\widehat{Y}_k(t)-Y_k(t)$ is equal to zero for each $t \in [0,T]$. This means that the model-assisted estimator $\widehat{\mu}_{MA}$ reduces in this case to the mean in the population  of the  predicted values $\widehat{Y}_k,$
\begin{align*}
\widehat{\mu}_{MA}(t) & =\frac{1}{N}\sum_{k\in U}\widehat{Y}_k(t), \quad t\in [0,T].
\end{align*}
Moreover, if only the intercept term is used, namely $Y_k(t)=\beta(t)+\varepsilon_{kt}$ for all $k\in U,$ then the estimator $\hat{\mu}_{MA}$ is simply  the well known H\'ajek estimator,
\begin{align*}
\widehat{\mu}_{MA}(t) &=\frac{\sum_{k\in s}\pi_k^{-1}Y_k(t)}{\sum_{k\in s}\pi_k^{-1}}, \quad t\in [0,T],
\end{align*}
which is sometimes preferred to the Horvitz-Thompson estimator (see e.g. \citet{SarndalLivre}, Chapter 5.7).
\end{Rem}

\begin{Rem}
 Estimator $\widehat{\mu}_{MA}(t)$ may also be obtained by using a calibration approach (\cite{DevilleSarndal1992}) which consists in looking for weights $w_{ks}, k\in s,$ that are as close as possible, according to some distance,  to the sampling weights $1/\pi_k$ while estimating exactly the population totals of the auxiliary information,
\begin{align*}
\sum_{k\in s}w_{ks}\mathbf{x}_k &=\sum_{k\in U}\mathbf{x}_k.
\end{align*}
Considering the chi-square distance leads to the following choice of weights 
\begin{align*}
w_{ks}=\frac{1}{\pi_k}-\left(\sum_{l\in s}\frac{\mathbf{x}_l}{\pi_l}-\sum_{l\in U}\mathbf{x}_l\right)'\left(\sum_{l\in s}\frac{\mathbf{x}_l\mathbf{x}_l'}{\pi_l}\right)^{-1}\frac{\mathbf{x}_k}{\pi_k}
\end{align*}
and the calibration estimator $\sum_sw_{ks}Y_k(t)/N$ for the mean $\mu(t)$ is equal  to $\widehat{\mu}_{MA}(t)$ defined in~(\ref{mu:ma}). 
\end{Rem}

\subsection{A regularized estimator for asymptotics}
The construction of the estimator $\widehat{\mu}_{MA}(t)$ is based on the assumption that  the matrix $\widehat{\mathbf{G}}$ is invertible. To show the uniform convergence, we  consider a modification of  $\widehat{\mathbf{G}}$ which will permit control of the expected norm of its inverse. Such a trick has already  been  used for example in  \citet{Bosq_2000} and \cite{Guillas2001}.
Since  $\widehat{\mathbf{G}}$ is a $p\times p$  symmetric and  non negative matrix it is possible to write it as follows
\begin{align*}
\widehat{\mathbf G} &=\sum_{j=1}^p\eta_{j,n} {\bf v}_{jn} {\bf v}'_{jn},
\end{align*} 
where $\eta_{j,n}$ is the  \textit{j}th eigenvalue,    $\eta_{1,n}\geq \cdots \geq \eta_{p,n}\geq 0$, and 
 ${\bf v}_{jn}$ is the corresponding orthonormal eigenvector.
Let us consider a real number $a>0$ and define the following regularized estimator of ${\bf G},$
\begin{align*}
\widehat{\mathbf G}_a &=\sum_{j=1}^p\max(\eta_{j,n},a) \  {\bf v}_{jn} {\bf v}'_{jn}.
\end{align*}
It is clear that $\widehat{\mathbf G}_a$ is always invertible and  
\begin{align}
\| \widehat{\mathbf G}_a^{-1} \|  & \leq a^{-1},
\label{borne_hatga}
\end{align}
where $\|.\|$ is the spectral norm for matrices.
Furthermore,  if $\eta_{p,n}\geq a$ then $\widehat{\mathbf G}=\widehat{\mathbf G}_a.$  
If $a>0$ is small enough, we show  under standard conditions on the moments of the variables $X_1, \ldots, X_p$ and on the first and second order inclusion probabilities that $\P(\widehat{\mathbf G}\neq\widehat{\mathbf G}_a)=\P(\eta_{p,n}< a) =O(n^{-1})$ (see Lemma \ref{norme Ga-G} in the Appendix).

Consequently,  it is possible to estimate the mean function $\mu_N(t)$ by the following estimator
 \begin{align}
\widehat{\mu}_{MA,a}(t)&=\frac{1}{N}\sum_{k\in U}\widehat{Y}_{k,a}(t)-\frac{1}{N}\sum_{k\in s}\frac{\widehat{Y}_{k,a}(t)-Y_k(t)}{\pi_k}, \quad t\in[0,T],
\label{mu:ma a}
 \end{align}
where $\widehat{Y}_{k,a}(t)=\mathbf{x}_k'\widehat{\boldsymbol\beta}_a(t)$ and $\widehat{\boldsymbol\beta}_a(t)=\widehat{\G}_a^{-1}\frac{1}{N}\sum_{k\in s}\frac{\mathbf{x}_kY_k(t)}{\pi_k}.$

\subsection{Discretized observations}

Note finally that with real data, we do not observe $Y_k(t)$ at all instants $t$ in $[0,T]$ but only for a finite set of $D$ measurement times, $0=t_1<...<t_D=T$. In functional data analysis, when the noise level is low and the grid of discretization points is fine, it is usual to perform a linear interpolation or to smooth  the discretized trajectories in order to obtain approximations of the trajectories at every instant $t \in [0,T]$ (see  \citet{Ramsay_Silverman_Livre}). 

If there are no measurement errors, if the trajectories are regular enough (but not necessarily differentiable) and if the grid of discretization points is dense enough,  \citet{cardot_josserand_2011} showed that  linear interpolation can provide sufficiently accurate approximations to the trajectories so that the approximation error can be neglected compared to the sampling error for the Horvitz-Thompson estimator. Note also that even if the observations are corrupted by noise, it has been shown by simulations  in \cite{CardotDegrasJosserand2012} that smoothing does really improve the accuracy of the Horvitz-Thompson estimator only when the noise level is high.

Thus, for each unit $k$ in  the sample $s$, we build the interpolated trajectory 
\begin{align*}
Y_{k,d}(t) &=Y_k(t_i)+\frac{Y_k(t_{i+1})-Y_k(t_i)}{t_{i+1}-t_i}(t-t_i) \quad t\in[t_i,t_{i+1}]
\end{align*}
and we define $\widehat{\boldsymbol\beta}_{a,d}(t)$ as the estimator of $\boldsymbol\beta(t)$ based on the discretized observations as follows
\begin{align*}
\widehat{\boldsymbol\beta}_{a,d}(t)&= \widehat{\bf G}_a^{-1}\frac{1}{N}\sum_{k\in s}\mathbf{x}_kY_{k,d}(t)\nonumber\\
&= \widehat{\boldsymbol\beta}_{a}(t_i)+\frac{\widehat{\boldsymbol\beta}_{a}(t_{i+1})-\widehat{\boldsymbol\beta}_{a}(t_i)}{t_{i+1}-t_i}(t-t_i).
\end{align*}
 
Therefore, the estimator of the mean population curve $\mu(t)$ based on the discretized observations is obtained by  linear interpolation  between  $\widehat{\mu}_{MA,a}(t_{i})$ and $\widehat{\mu}_{MA,a}(t_{i+1})$.  For $t \in [t_i,t_{i+1}]$, 
\begin{align}
\widehat{\mu}_{MA,d}(t)
& =  \frac{1}{N}\sum_{k\in U}\widehat{Y}_{k,d}(t)-\frac{1}{N}\sum_{k\in s}\frac{(\widehat{Y}_{k,d}(t)-Y_{k,d}(t))}{\pi_k} \nonumber \\
 &=  \widehat{\mu}_{MA,a}(t_i)+\frac{\widehat{\mu}_{MA,a}(t_{i+1})-\widehat{\mu}_{MA,a}(t_i)}{t_{i+1}-t_i}(t-t_i) 
\label{diffmuad}
\end{align}
 where $\widehat{Y}_{k,d}(t)=\mathbf{x}'_k\widehat{\boldsymbol\beta}_{a,d}(t). $
 \section{Asymptotic properties under the sampling design}
 All the proofs are postponed in an Appendix. 
 \subsection{Assumptions}
 To derive the asymptotic properties under the sampling design $p(\cdot)$ of  $\widehat{\mu}_{MA,d}$ we must suppose that  both the sample size and the population size become large. More precisely, we consider the superpopulation framework introduced by \citet{isaky_fuller_1982} with  a sequence of growing and nested populations $U_N$ with size $N$ tending to infinity and a sequence of samples $s_N$ of size $n_N$ drawn from $U_N$ according to the sampling design $p_N(s_N)$. The first and second order inclusion propabilities are respectively denoted by $\pi_{kN}$ and $\pi_{klN}$. For simplicity of notations and when there is no ambiguity, we drop the subscript $N$. To prove our asymptotic results we need the following assumptions.
\begin{itemize}
\item[\bf A1.] We assume that $\displaystyle\lim_{N\rightarrow \infty} \frac{n}{N}=\pi \in (0,1).$
\item[\bf A2.] We assume that $\displaystyle\min_{k \in U} \pi_k\geq\lambda>0$, $\displaystyle\min_{k \neq l} \pi_{kl}\geq\lambda^*>0$ and

 $\displaystyle\limsup_{N \rightarrow \infty}n \max_{k\neq l\in U}\vert \pi_{kl}-\pi_k\pi_l\vert<C_1<\infty$
 
\item[\bf A3.] There are two positive constants  $C_2$ and $C_3$ and $1 \geq \beta>1/2$ such that, for all $N$ and for all $(r,t)\in[0,T]\times [0,T]$,

$$\frac{1}{N}\sum_{k\in U}Y_k(0)^2<C_2 \quad \mbox{and} \quad \frac{1}{N}\sum_{k\in U}\{Y_k(t)-Y_k(r)\}^2<C_3\vert t- r\vert^{2\beta}.$$ 

\item[\bf A4.] We assume that there is a positive constant  $C_4$ such that for all $k\in U,$  $\Vert\mathbf{x}_k\Vert^2<C_4. $
\item[\bf A5.] We assume that, for $N>N_0,$ the matrix  $\mathbf{G}$ is invertible and  that the number $a$ chosen before satisfies  $\| \mathbf{G}^{-1}\| < a^{-1}.$
\end{itemize}

Assumptions {\bf A1} and {\bf A2} are classical hypotheses in survey sampling and deal with the first and second order inclusion probabilities. They are satisfied for many usual sampling designs with fixed size (see for example \cite{Hajek1981}, \cite{RobinsonSarndal83} and \cite{breidt_opsomer_2000}).

Assumption {\bf A3} is a minimal regularity condition already required in \cite{cardot_josserand_2011}. Even if pointwise consistency, for each fixed value of $t$, can be proved without any condition on the H\"older coefficient $\beta$,  this regularity condition is necessary to get a uniform convergence result. A counterexample is given in \cite{Hahn77} when $\beta \leq  1/2.$ More precisely it is shown that the sample mean i.i.d copies of a uniformly bounded continuous random function defined on a compact interval may not satisfy the Central Limit Theorem in the space of continuous functions. The hypothesis $\beta>1/2$ also implies that the trajectories of the residual processes $\epsilon_{kt},$ see (\ref{model:xi}), are regular enough (but not necessarily differentiable).  Assumption {\bf A4} could certainly be weakened at the expense of longer proofs.
Assumption {\bf A5} means that for all $\mathbf{u} \in \mR,$ with $\mathbf{u} \neq 0,$ we have $\mathbf{u}'\mathbf{Gu} \geq a \mathbf{u}'\mathbf{u}.$ The same kind of assumption  is required in \cite{isaky_fuller_1982} to get the pointwise convergence in probability whereas 
\cite{RobinsonSarndal83} introduce a much stronger condition (condition A7 in their article) which directly deals with the mean square convergence of the estimator of the vector $\boldsymbol{\beta}$ of regression coefficients. 

  \subsection{Uniform consistency of $\hat{\mu}_{MA,d}$}
We aim at showing  that $\widehat{\mu}_{MA,d}$ is uniformly consistent for $\mu,$ namely  that, for all $\varepsilon>0$,  
\begin{align*}
\P \left(\sup_{t\in [0,T]}|\widehat{\mu}_{MA,d}(t)-\mu(t)|>\varepsilon \right) &\rightarrow 0,
\end{align*}
when $N$ tends to infinity.
The suitable space for proving the uniform convergence is the space of continuous functions on $[0,T]$, denoted by $C[0,T]$, equipped with its natural distance $\rho$;   for two elements  $f,g \in C[0,T],$  the distance between $f$ and $g$ is  $\rho(f,g)=\sup_{t\in [0,T]}|f(t)-g(t)|.$
It results that the uniform consistency of $\widehat{\mu}_{MA,d}$ is simply the convergence in probability of $\widehat{\mu}_{MA,d}$ to $\mu$ in the space $C[0,T].$  Remark that with assumption {\bf A3}  the trajectories $Y_k$ are continuous for all $k\in U,$ and thus the mean curve $\mu$ belongs to $C[0,T]$ as well as its estimator  $\widehat{\mu}_{MA,d},$ by construction.  

We first state the uniform consistency of the estimator $\widehat{\boldsymbol{\beta}}_{a,d}(t)$ towards its population counterpart $\boldsymbol{\tilde\beta}(t)$ under conditions on the number and the repartition of discretization points.
\begin{Prop}
Let assumptions (A1)-(A5) hold. If the discretization scheme satisfies 
$\max_{i \in \{1,..,D_N-1\}}\vert t_{i+1}-t_i\vert^{2\beta}=o(n^{-1})$ then there is a constant $C>0$ such that, for all n,
\begin{align*}
\sqrt{n} \ \E_p\left\lbrace \sup_{t\in [0,T]} \left\| \widehat{\boldsymbol{\beta}}_{a,d}(t)-\boldsymbol{\tilde\beta} (t)\right\|  \right\rbrace  &\leq  C.
\end{align*}
\label{convergencebeta}
\end{Prop}

\noindent We can now state a similar type of result for the estimator of the mean function.
\begin{Prop}
Let assumptions (A1)-(A5) hold. If the discretization scheme satisfies 
$\max_{i \in \{1,..,D_N-1\}}\vert t_{i+1}-t_i\vert^{2\beta}=o(n^{-1})$ then there is a constant $C>0$ such that, for all n,
\begin{align*}
\sqrt{n} \ \E_p\left\lbrace \sup_{t\in [0,T]}\mid\widehat{\mu}_{MA,d}(t)-\mu(t)\mid\right\rbrace &\leq  C.
\end{align*}
\label{convergence mu}
\end{Prop}

We deduce from  Proposition \ref{convergence mu} that  estimator $\widehat{\mu}_{MA,d}(t)$ is asymptotically unbiased as well as  design consistent.
Note that the approximation error (with linear interpolation) is negligible, compared to the sampling variability, under the additional assumption on the repartition of  the discretization points.
This assumption also tolds us that less discretization points are required for smoother trajectories.

Let us also remark that, for each $t,$
\begin{align}
\widehat{\mu}_{MA,a}(t)-\widetilde{\mu}(t)&=  \frac{1}{N} \sum_{k \in U} \left(1-\frac{\1_k}{\pi_k}  \right) \mathbf{x}_k' \left( \widehat{\boldsymbol{\beta}}_{a}(t) -\widetilde{\boldsymbol{\beta}}(t)\right),
\label{def:mutilde-muhat}
\end{align}
where $\1_k$ is the sample membership, so that it is  not difficult to prove, under previous assumptions and by using lemma \ref{espbetatinstantt} in the Appendix, that for all $t \in [0,T]$,
\begin{align}
\sqrt{n}\left(\widehat{\mu}_{MA,d}(t)-\widetilde{\mu}(t)\right) &=o_p(1).
\label{approx_hatmu_tildemu}
\end{align} 

\subsection{Covariance function estimation under the sampling design}
We undertake in this section a detailed study of the  covariance function  of estimator $\hat{\mu}_{MA,d}.$ The covariance function is computed with respect to the sampling design $p(\cdot)$ and from relation (\ref{diffmuad}), we can deduce that $\widehat{\mu}_{MA,d}$ is a nonlinear function of Horvitz-Thompson estimators, so the usual Horvitz-Thompson covariance formula given by (\ref{covar_ht}) can not be used anymore. Nevertheless, in light of relation (\ref{approx_hatmu_tildemu}), the covariance function of $\widehat{\mu}_{MA,d}$ between two instants $r$ and $t$ may be approximated by the covariance $\mbox{Cov}_p(\tilde{\mu}(r), \tilde{\mu}(t)),$  which in turn is equal to the Horvitz-Thompson covariance applied to the residuals $Y_k-\tilde{Y}_k.$ Let us denote by $\gamma_{MA}$  the approximative covariance function of $\widehat{\mu}_{MA,d}$  defined as follows 

\begin{align}
 \gamma_{\textrm{MA}}(r,t)  & =  \frac{1}{N^2}\mbox{Cov}_p\left(\sum_{k\in s} \frac{Y_k(r)-\tilde{Y}_k(r)}{\pi_k},\sum_{k\in s} \frac{Y_k(t)-\tilde{Y}_k(t)}{\pi_k}\right)\nonumber\\
  &=  \frac{1}{N^2}\sum_{k\in U}\sum_{l\in U}(\pi_{kl}-\pi_k \pi_l)\frac{Y_k(r)-\tilde{Y}_k(r) }{\pi_k}\frac{ Y_l(t)-\tilde{Y}_l(t) }{\pi_l}, \quad r,t \in [0,T].
   \label{var_ma_asymp}
\end{align}

This approximation explains that  model-assisted estimators will perform much better than Horvitz-Thompson estimators if the  residuals $Y_k(t)-\tilde{Y}_k(t)$ are small compared to $Y_k(t)$. The covariance function $\gamma_{\textrm{MA}}(r,t) $ can be estimated  by the Horvitz-Thompson variance estimator for the estimated residuals $Y_{k,d}(t)-\widehat{Y}_{k,d}(t),$ 
 \begin{align}
    \widehat{\gamma}_{\textrm{MA},d}(r,t) &= \frac{1}{N^2} \sum_{k,l \in s}\frac{\pi_{kl}-\pi_k \pi_l}{ \pi_{kl}} \cdot\frac{ Y_{k,d}(r)-\widehat{Y}_{k,d}(r) }{\pi_k}\cdot\frac{  Y_{l,d}(t)-\widehat{Y}_{l,d}(t) }{\pi_l}, \quad r,t \in [0,T],
     \label{cov MA}    
  \end{align}
 where $ \widehat{Y}_{k,d}(t) =  {\bf x}_k'\widehat{\boldsymbol \beta}_{a,d}(t).$\\
To prove the consistency of the covariance estimator $ \widehat{\gamma}_{\textrm{MA},d}(r,t)$, let us introduce additional assumptions that involve  higher-order inclusion probabilities as well as conditions on the  fourth order moments of the trajectories.

\begin{itemize}
\item[{\bf A6.}] We assume that  
\begin{align*}
\lim_{N\rightarrow\infty}\max_{(k,l,k',l')\in D_{4,n}}\vert \E_p\{(\1_{kl}-\pi_{kl})(\1_{k'l'}-\pi_{k'l'})\}\vert&=0
\end{align*}
 where $D_{t,N}$ is the set of all distinct $t$-tuples $(i_1,...,i_t)$ from $U_N$ and $\1_{kl}=\1_{k}\1_{l}.$

\item[{\bf A7.}] There are two positive constants $C_5$ and $C_6$ such that $N^{-1}\sum_UY_k(0)^4<C_5$ and $N^{-1}\sum\{Y_k(t)-Y_k(r)\}^4<C_6\vert t-r\vert^{4\beta}$, for all $(r,t)\in[0,T]^2$
\end{itemize}

Condition {\bf A6} has already been assumed by \cite{breidt_opsomer_2000}  in a nonparametric model-assisted context and in \citet{cardot_josserand_2011} for Horvitz-Thompson estimators. It can be checked that it is fulfilled for simple random sampling without replacement (SRSWOR) or stratified sampling with SRSWOR within each strata. More generally, it is fulfilled for high entropy sampling designs. \cite{BLRG2012} prove that it is fulfilled for the  rejective sampling whereas  \cite{CGL2012} check that it is true for  sampling designs, such as Sampford sampling, whose Kullback-Leibler divergence with respect to rejective sampling, tends to zero when the population size increases.

\begin{Prop} 
Assume (A1)-(A7) hold and the sequence of discretization schemes satisfy ${\lim_{N\rightarrow\infty}\max_{i \in \{1,..,D_N-1\}}\vert t_{i+1}-t_i\vert=0}.$ Then, as $N$ tends to infinity, we have for all $(r,t) \in [0,T]^2,$
\begin{eqnarray*}
n \  \E_p \left\{ \left| \widehat{\gamma}_{\textrm{MA},d}(r,t)-\gamma_{\textrm{MA}}(r,t) \right| \right\} &\rightarrow& 0 
\end{eqnarray*}
and 
\begin{eqnarray*}
n \ \E_p \left\{ \sup_{t \in [0,T]} \left| \widehat{\gamma}_{\textrm{MA},d}(t,t)-\gamma_{\textrm{MA}}(t,t) \right| \right\} &\rightarrow& 0. 
\end{eqnarray*}
\label{prop:gammaconv} 
\end{Prop}

Since $n \gamma_{\textrm{MA}}(r,t)$ remains bounded, the previous proposition  tells us that $\widehat{\gamma}_{\textrm{MA},d}$ is  consistent pointwise and the variance function estimator is uniformly convergent.  Note also that the condition on the number of discretization points is much weaker than in Proposition~\ref{convergence mu} because we do not give here rates of convergence. To obtain such rates, we would also need  to impose additional assumptions on the sampling design.
\subsection{Asymptotic normality and confidence bands}
We assume  a supplementary assumption  in order to get the asymptotic normality of the functional estimator $\hat{\mu}_{\textrm{MA},d}$ in the space of continuous functions.
\begin{itemize}
\item[{\bf A8.}] We assume that for each fixed value of $t \in [0,1]$, 
\begin{align*}
\{\gamma_{\textrm{MA}}(t,t)\}^{-1/2}\left(\tilde{\mu}(t)-\mu(t) \right) &\rightarrow   {\cal N}(0,1)
\end{align*}
 in distribution when $N$ tends to infinity. 
\end{itemize}
This assumption is satisfied for usual sampling designs (see {\it e.g.} \cite{FullerLivre}, Chapter 2.2). Note that using relation (\ref{approx_hatmu_tildemu}), we can write for any fixed value $t \in [0,T],$
\begin{align*}
\widehat{\mu}_{\textrm{MA},d}(t)-\mu(t) &=\tilde{\mu}(t)-\mu(t)+o_p(n^{-1/2}),
\end{align*}
and deduce that  $\sqrt{n}\left(\widehat{\mu}_{\textrm{MA},d}(t)-\mu(t) \right)$ is also pointwise asymptotically Gaussian when conditions of Proposition \ref{convergencebeta} hold. Let us  state now a much stronger result which indicates that the convergence to a Gaussian distribution also occurs for the trajectories, in the space of continuous functions (see \cite{Billingsley1968}, Chapter 2). 

\begin{Prop}
Assume (A1)-(A5) and (A8) hold. If the discretization scheme satisfies 
$\max_{i=\{1,..,D_N-1\}}\vert t_{i+1}-t_i\vert^{2\beta}=o(n^{-1}),$ we have when $n$ tends  to infinity
\begin{align*}
\sqrt{n}\left\lbrace \widehat{\mu}_{\textrm{MA},d}-\mu \right\rbrace &\rightsquigarrow Z
\end{align*} 
\textit{where $\rightsquigarrow$ indicates the convergence in distribution in $C[0,T]$ with the uniform topology and
$Z$ is a Gaussian process taking values in $C[0,T]$ with mean 0 and covariance function 
$\gamma_Z(r,t) = \lim_{n \rightarrow + \infty} n \gamma_{\textrm{MA}}(r,t).$}
\label{prop:CLTdansC}
\end{Prop}

The "sup" functional defined on the space of continuous functions being continuous, the Proposition~\ref{prop:CLTdansC} implies that the real random variable 
$\sup_t |\sqrt{n}\left\lbrace \widehat{\mu}_{\textrm{MA},d}(t)-\mu(t) \right\rbrace|$ converges in distribution to $\sup_t |Z(t)|$.
We thus consider confidence bands  for $\mu$ of the form 
\begin{align}\label{SCB}
\left\{ \left[ \widehat{\mu}_{\textrm{MA},d}(t) \pm c \, \frac{  \widehat{\sigma}(t)}{\sqrt{n}} \right], \: t\in [0,T] \right\} 
\end{align}
where $c$ is a suitable number and $\widehat{\sigma}(t) = \sqrt{n \widehat{\gamma}_{\textrm{MA},d}(t,t)}$. Note that the fact that $\mu$ belongs to the confidence band defined in (\ref{SCB}) is equivalent to 
\begin{align*}
\sup_{t \in [0,T]} \frac{\sqrt{n}}{\widehat{\sigma}(t)} \left| \widehat{\mu}_{\textrm{MA},d}(t) - \mu(t) \right| &\leq c.
\end{align*}

Given a confidence level $1-\alpha \in (0,1)$, one way to build such confidence band, that is to say one way  to find an adequate value for $c_\alpha,$ is to perform simulations of a centered Gaussian functions $\widehat{Z}$ defined on $[0,T]$ with mean 0 and covariance function $n \widehat{\gamma}_{\textrm{MA},d}(r,t)$ and then compute the quantile of order $1-\alpha$ of $\sup_{t \in [0,T]} \left| \widehat{Z}(t)/\widehat{\sigma}(t) \right|.$ In other words, we look for a constant $c_\alpha$, which is in fact a random variable since it depends on the estimated covariance function $\widehat{\gamma}_{\textrm{MA},d},$ such that
\begin{align*}
\P \left(  |  \widehat{Z}(t) | \leq c_\alpha  \frac{  \widehat{\sigma}(t)}{\sqrt{n}}, \ \forall t \in [0,T] \mid \widehat{\gamma}_{\textrm{MA},d} \right) &= 1 - \alpha
\end{align*}
The asymptotic coverage of this simulation based procedure has been rigorously studied for the Horvitz-Thompson estimators of the mean of sampled and noisy trajectories in  \cite{CardotDegrasJosserand2012} whereas \cite{CDGJL2012} have successfully employed this approach on real load curves with model-assisted estimators.
The next proposition, which can be seen as a functional version of Slutsky's Lemma, provides a rigorous justification of this latter technique.

\begin{Prop}\label{conditional weak functional convergence for Gaussian processes}
Assume (A1)-(A8) hold and the discretization scheme satisfies 
$\max_{i \in \{1,..,D_N-1\}}\vert t_{i+1}-t_i\vert^{2\beta}=o(n^{-1}).$ 

Let  $Z$ be a Gaussian process with mean zero and covariance function $\gamma_Z$ (as  in Proposition \ref{prop:CLTdansC}). 
Let $(\widehat{Z}_N)$ be a sequence of processes such that for each $N$, conditionally on the estimator $\widehat{\gamma}_{\textrm{MA},d}$ defined in (\ref{cov MA}), $\widehat{Z}_N$ is Gaussian with mean zero and covariance $n \widehat{\gamma}_{\textrm{MA},d}$. Suppose that $\gamma_Z(t,t)$ is a continuous function and $\inf_{t} \gamma_Z(t,t) > 0.$
Then, as $N\to\infty$, the following convergence holds uniformly in $c$, 
\begin{align*}
 \mathbb{P} \left(  |\widehat{Z}_N(t)| \le c\, \widehat{\sigma} (t) ,  \: \forall t \in [0,T] \, \big| \,  \widehat{\gamma}_{\textrm{MA},d} \right) &\to   \mathbb{P} \left(  |Z(t)| \le c \, \sigma (t) , \: \forall t \in [0,T] \right),
 \end{align*}
 where $\widehat{\sigma}(t) = \sqrt{n \widehat{\gamma}_{\textrm{MA},d}(t,t)}$ and $\sigma(t) = \sqrt{\gamma_Z(t,t)}.$
\end{Prop}
As in \cite{CardotDegrasJosserand2012}, it is possible to deduce from the previous proposition that the chosen value $\widehat{c}_\alpha = c_\alpha(\widehat{\gamma}_{\textrm{MA},d})$ provides asymptotically the desired coverage since it satisfies 
 \begin{align*}
\lim_{N\to\infty} \mathbb{P} \left( \mu(t)\in \left[ \widehat{\mu}_{\textrm{MA},d}(t) \pm \widehat{c}_\alpha \frac{\widehat{\sigma}(t)}{\sqrt{n}}  \right], \ \forall t\in [0,T] \right) &= 1-\alpha.
\end{align*}

\section{An illustration on electricity consumption curves}
We consider, as in \cite{CDGJL2012}, a population of $N=15069$ electricity consumption curves, measured every 30 minutes over a period of one week. Each element $k$ of the population is thus a vector with size 336, denoted by $(Y_k(t), t \in \{1, \ldots, 336\}).$ The auxiliary information $X$ of values $x_k, k \in U$ is simply the mean consumption of each meter $k\in U$ recorded during the  week before the sample is drawn. As shown in Figure~\ref{fig:tempcorr}, the real variable $X$ is strongly correlated with the consumption at each instant $t$ of the current period of estimation so that a linear model with a functional response is well adapted for model-assisted estimation.

We draw samples $s_i$ of size $n,$ for $i=1, \ldots, I=10000$ with simple random sampling without replacement (SRSWOR) so that $\pi_k = n/N$ for $k \in \{1, \ldots, N\}$. We define, for each sample $s_i,$ the model-assisted estimator of the mean curve, 
\begin{align}
\widehat{\mu}_{MA,d}^{(i)}(t) &=\frac{1}{N}\sum_{k\in U}\widehat{Y}^{(i)}_k(t)-\frac{1}{N}\sum_{k\in s_i}\frac{\widehat{Y}^{(i)}_k(t)-Y_k(t)}{n/N}
\label{def:MAestsimu}
\end{align}
where   $\mathbf{x}'_k = (1,x_k)$, $\widehat{Y}^{(i)}_k(t)=\mathbf{x}_k'\widehat{\boldsymbol\beta}^{(i)}(t),$ and $\widehat{\boldsymbol\beta}^{(i)}(t)=\widehat{\mathbf{G}}^{-1}\frac{1}{N}\sum_{k \in s_i}\frac{\mathbf{x}_kY_k(t)}{n/N}$ for $t \in \{1, \ldots, 336\}.$ \cite{CDGJL2012} noted that, for the same sample size, the mean square error of estimation of the mean curve is divided by four compared to the Horvitz-Thompson estimator with SRSWOR when considering the model-assisted estimator defined in~(\ref{def:MAestsimu}).  There is only one covariate in this study and we did not encounter any problem with the invertibility of matrix $\widehat{\mathbf{G}},$ the value of parameter $a$ is thus  $a=0$.

\begin{figure}
\begin{center}
\includegraphics[width=\textwidth]{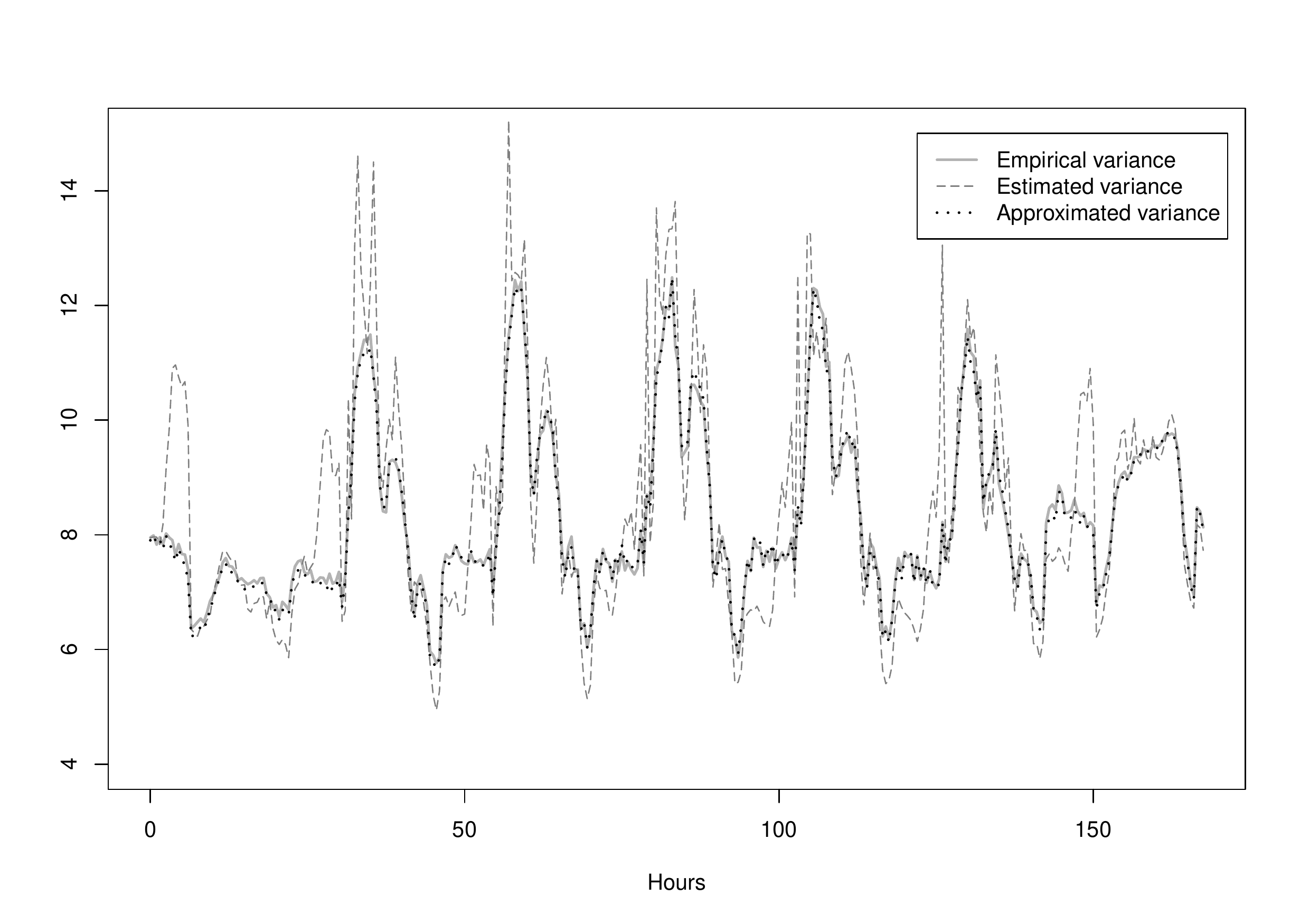}
\caption{Empirical variance function $\gamma_{emp}$, approximated variance $\gamma_{MA}$ and estimated variance $\hat{\gamma}_{MA,d}$ obtained with a sample of size $n=1500$.}  
\label{fig:Estvariance}
\end{center}
\end{figure}

We also define $\hat{\overline{\mu}}(t)=\frac{1}{I}\sum_{i=1}^I\hat{\mu}_{MA,d}^{(i)}(t),$ $t \in \{1, \ldots, 336\}.$
The true variance function of the model-assisted estimator being unknown, we approximate it with a Monte Carlo approach based on the $I=10000$ samples drawn with simple random sampling without replacement. The approximation to the true variance function is thus given by
\begin{align}
\gamma_{emp}(r,t) &=\frac{1}{I}\sum_{i=1}^I(\hat{\mu}_{MA,d}^{(i)}(t)-\hat{\overline{\mu}}(t))(\hat{\mu}_{MA,d}^{(i)}(r)-\hat{\overline{\mu}}(r))
\end{align}
for $(r,t) \in \{1, \ldots, 336\}$.

The following quadratic loss criterion which measures a relative error is used to evaluate, for each sample, the accuracy of the variance estimator defined in (\ref{cov MA}),
\begin{align}
E_r(\hat{\gamma}_{MA,d}) &=\frac{1}{336} \sum_{t=1}^{336} \frac{\vert \hat{\gamma}_{MA,d}(t,t)-\gamma_{emp}(t,t)\vert^2}{\gamma_{emp}(t,t)^2 }dt \quad
\label{def:L2loss}
\end{align}
We also decompose, over the $I=10000$ estimations, the relative mean square error (RMSE) of the estimator into an approximation error ($RB(\hat{\gamma}_{MA,d})^2$) and a variance term ($VR(\hat{\gamma}_{MA,d})$) that can be related to the sampling error,
\begin{align*}
RMSE(\hat{\gamma}_{MA,d})&= \frac{1}{I} \sum_{i=1}^IE_r^{(i)}(\hat{\gamma}_{MA,d})\nonumber\\
&= RB(\hat{\gamma}_{MA,d})^2+VR(\hat{\gamma}_{MA,d})
\end{align*}
where $E_r^{(i)}(\hat{\gamma}_{MA,d})$ is the value of $E_r(\hat{\gamma}_{MA,d})$ for the $i$th sample. The relative bias of the estimator $\hat{\gamma}_{MA,d}$ may be written as
\begin{align*}
RB(\hat{\gamma}_{MA,d})^2&=\frac{1}{336}\sum_{t=1}^{336}\left(\frac{\overline{\hat{\gamma}}_{MA,d}(t,t)-\gamma_{emp}(t,t)}{\gamma_{emp}(t,t)}\right)^2
\end{align*} 
where $\overline{\hat{\gamma}}_{MA,d}(t,t)=\frac{1}{I}\sum_{i=1}^I\hat{\gamma}_{MA,d}^{(i)}(t,t)$.

\begin{table}[htbp]
\begin{center}
\begin{tabular}{|c|c|c|c|c|c|c|c|}
   \hline
   Sample size &$RMSE(\hat{\gamma}_{MA,d})$&$RB(\hat{\gamma}_{MA,d})^2$&\multicolumn{5}{|c|}{$E_r(\hat{\gamma}_{MA,d})$}\\
   \cline{4-8}
    && & $q_5$&$q_{25}$& Median & $q_{75}$ & $q_{95}$ \\      
     \hline
   250&0.1315&0.0027&0.0264&0.0455&0.0707&0.117&0.4945\\
     \hline
   500&0.0697&0.0016&0.0166&0.029&0.0459&0.0794&0.1945\\
 \hline
   1500&0.0238&0.0003&0.0076&0.0125&0.0186&0.028&0.0569\\
 \hline

\end{tabular}
\end{center}
\caption{Summary statistics for $E_r(\hat{\gamma}_{MA,d},\gamma_{emp}),$ with  I=10000 samples.}
\label{largeurbc}
\end{table}

The RMSE as well as the approximation error and statistics (quantiles) for $E_r$ are given in Table~\ref{largeurbc}. We can note that logically the RMSE decreases as the sample size increases and that even for moderate sample sizes, the estimations are rather precise. A closer look on how the RMSE is decomposed reveals that estimation error is mainly due to the sampling error, via the variance term whereas the approximation error term $RB(\hat{\gamma}_{MA,d})^2$ is negligible.  This fact can be observed  in Figure~\ref{fig:Estvariance} were we plot  the true variance function $\gamma_{emp}$ over the considered period, its approximation $\gamma_{MA}$ as well as an estimation $\widehat{\gamma}_{MA,d}$ with a sample with size $n=1500$, whose error according to criterion (\ref{def:L2loss}) is close to the mean error ($E_r \approx 0.02$).

\begin{figure}
\begin{center}
\includegraphics[width=\textwidth]{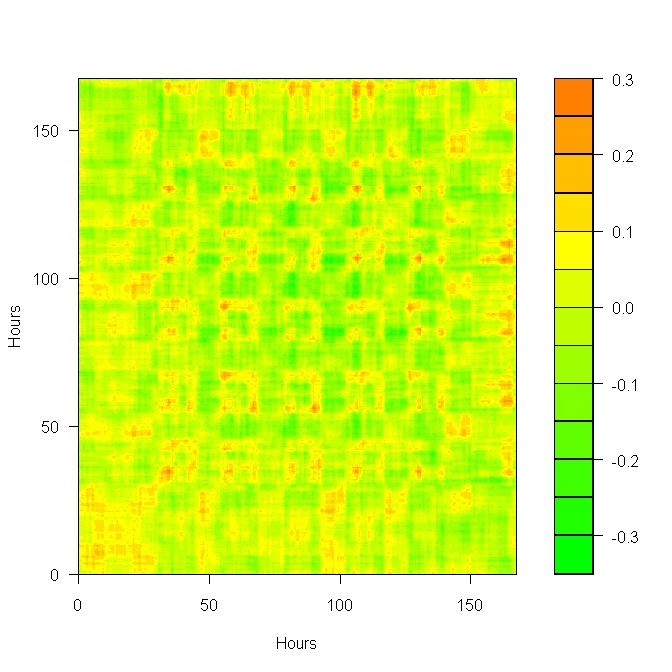}
\caption{(Approximation error) difference between the covariance function and its approximation, $\gamma_{emp}(t,r)-{\gamma}_{MA}(t,r),$ for a sample with size  $n=1500$.}  
\label{fig:erreurcovariance}
\end{center}
\end{figure}

\begin{figure}
\begin{center}
\includegraphics[width=\textwidth]{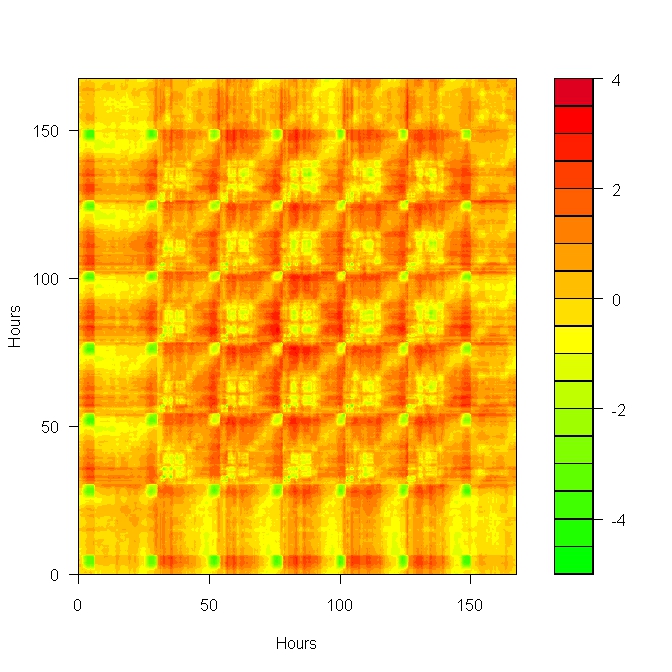}
\caption{(Sampling error) difference between the approximated covariance function and its estimation, ${\gamma}_{MA}(t,r)-\hat{\gamma}_{MA,d}(t,r),$ for a sample with size $n=1500$.}  
\label{fig:erreurcovarianceapproximee}
\end{center}
\end{figure}

We have  also plotted in Figure~\ref{fig:erreurcovariance} the difference between the empirical covariance function $\gamma_{emp}$ and its approximation $\gamma_{MA}$ and  in Figure~\ref{fig:erreurcovarianceapproximee} the difference between  $\gamma_{MA}$ and its estimation $\widehat{\gamma}_{MA,d}$ for a sample with size $n=1500$ whose error, $E_r \approx 0.02,$ is close to the mean value. Once again, it is clearly seen that the approximation error to the true covariance function (see Figure~\ref{fig:erreurcovariance}) is much smaller than the sampling error (see Figure~\ref{fig:erreurcovarianceapproximee}). We can also remark some strong periodic pattern which reflects  the natural daily periodicity in the electricity consumption behavior and that is related to the temporal correlation of the unknown residual process $\epsilon_{kt}$ defined in (\ref{model:xi}).

\section{Concluding remarks}

We have made in this paper an asymptotic study of model-assisted estimators, with linear regression models with functional response,  when the target is the mean of functional data with discrete observations in time.
This work can be extended in many directions. 
For example, one could consider more sophisticated regression models than model (\ref{model:xi}) such as non linear or nonparametric models with functional response by adapting, in a survey sampling context, models studied in the functional data analysis literature (see \cite{CMW2004}, \cite{Cardot2007}, or  \cite{FLTV2011}). However, one important drawback of such more sophisticated approaches is that they would require to  know $\mathbf{x}_k$ for all the units $k$ in the population as opposed to only their population totals.

An interesting direction for future investigation would be to consider noisy and possibly sparse measurements in time. For the Horvitz-Thompson estimator, local polynomials are employed in  \cite{CardotDegrasJosserand2012} in order to first smooth the trajectories and it would certainly be possible to adapt the techniques developed in this work to the model-assisted estimation procedure.

Another promising direction for future research would be to adapt model-assisted estimators for  time-varying samples. When one works with large networks of sensors it can be possible to consider a sequence of samples $s(t)$ that evolve along time. A preliminary work (see \cite{Degras2012}), which focuses  on Horvitz-Thompson estimators and stratified sampling clearly shows that such  time-varying samples can  outperform sampling designs that are fixed in time. 

\bigskip

\noindent{\bf Acknowledgements.} We thank the two anonymous referees for a careful reading and interesting suggestions that have permitted a great  improvement of  the original manuscript.

\appendix
\section{Proofs}
Throughout the proofs we use the letter $C$ to denote a generic constant whose value may vary from place to place. We  also denote by $\alpha_k=\frac{\1_k}{\pi_k}-1,$ $k\in U$ and by $\Delta_{kl}=\pi_{kl}-\pi_k\pi_l, $ $k,l\in U.$

\subsection{Some useful Lemmas}
Note that the result showed in the first Lemma  is sometimes stated as an assumption (see \textit{e.g} \cite{RobinsonSarndal83}). It is used to prove  the convergence of the estimator of the mean  in terms of mean square error. 
\begin{Lem}
Let assumptions  (A1), (A2) and (A4), (A5) hold. Then, there is a constant $C$ such that
\begin{align*}
n\ \E_p \left(\Vert\widehat{{\bf G}}_a^{-1}-{\bf G}^{-1}\Vert^2\right) &\leq C.
\end{align*}
\label{norme Ga-G}
\end{Lem}


\begin{Proof}. The proof follows the lines of (\cite{Bosq_2000}, Theorem 8.4) and (\cite{cardot_chaouch_goga_labruere_2010}, Proposition 3.1). Using assumption (A5) and inequality  (\ref{borne_hatga}), we have
\begin{align*}
\Vert \widehat{{\bf G}}_a^{-1}-{\bf G}^{-1} \Vert&\leq  \Vert \widehat{\bf G}_a^{-1}\Vert.\Vert\widehat{\bf G}_a-{\bf G}\Vert .\Vert {\bf G}^{-1}\Vert\nonumber\\
&\leq a^{-2}\Vert\widehat{\bf G}_a-{\bf G}\Vert, 
\end{align*}
which implies
\begin{align}
\E_p\left(\Vert\widehat{\bf G}_a^{-1}-{\bf G}^{-1}\Vert^2\right)&\leq a^{-4}\E_p\left(\Vert\widehat{\bf G}_a-{\bf G}\Vert^2\right). 
\label{equation 1}
\end{align}
To bound $\E_p\left(\Vert\widehat{\bf G}_a-{\bf G}\Vert^2\right),$ we use the following decomposition. 
\begin{align}
 \E_p\left(\Vert\widehat{\G}_a-\G\Vert^2\right)
&= \E_p\left(\Vert\widehat{\G}_a-\G\Vert^2\1_{\{\widehat{\G}_a=\widehat{\G}\}}\right)+\E_p\left(\Vert\widehat{\G}_a-\G\Vert^2\1_{\{\widehat{\G}_a\neq\widehat{\G}\}}\right) \nonumber\\
&\leq \E_p\left(\Vert\widehat{\G}-\G\Vert^2\right)+2\E_p\left(\Vert\widehat{\G}_a-\widehat{\G}\Vert^2\1_{\{\widehat{\G}_a\neq\widehat{\G}\}}\right)+ 2\E_p\left(\Vert\widehat{\G}-\G\Vert^2\1_{\{\widehat{\G}_a\neq\widehat{\G}\}}\right)\nonumber\\
& \leq   3\E_p(\Vert\widehat{\G}-\G\Vert^2) +2\E_p\left(\Vert\widehat{\G}_a-\widehat{\G}\Vert^2\1_{\{\widehat{\G}_a\neq\widehat{\G}\}}\right).
\label{esp_ga_g}
\end{align}
To bound the first term from the right-side of (\ref{esp_ga_g}), we use the fact that the spectral norm is majored by the trace norm  $||\cdot||_2$ defined by  $\left\| {\bf A} \right\|_2^2=tr({\bf A}'{\bf A}).$ Next, we show (see also \cite{cardot_chaouch_goga_labruere_2010}, proof of Proposition 3.1) that,
\begin{align}
\E_p\Vert \widehat{\G}-\G\Vert_2^2&=O(n^{-1}). \label{equation 2}
\end{align}
 We have, with assumptions (A1), (A2) and (A4) that, 
\begin{align*}
\E_p\Vert \widehat{\G}-\G\Vert_2^2&= \frac{1}{N^2}\E_p\left(\sum_{k\in U}\sum_{l\in U}\alpha_k\alpha_l\text{tr}[\mathbf{x}_k\mathbf{x}_k'\mathbf{x}_l\mathbf{x}_l']\right)\nonumber\\
&\leq \frac{1}{N^2}\frac{1}{\lambda}\sum_{k\in U}\Vert\mathbf{x}_k\mathbf{x}_k'\Vert_2^2+\max_{k\neq l\in U}\vert\Delta_{kl}\vert\frac{1}{N^2\lambda^2}\sum_{k\in U}\sum_{l\in U} \Vert\mathbf{x}_k\Vert^2\Vert\mathbf{x}_l\Vert^2\nonumber\\
&\leq \frac{1}{n}\left(\frac{n}{N}\frac{1}{\lambda}+n\max_{k\neq l\in U}\vert\Delta_{kl}\vert\frac{1}{\lambda^2}\right)C^2_2\nonumber\\
&\leq  \frac{C}{n}.
\end{align*}
On the other hand,
$$
\E_p\left(\Vert\widehat{\G}_a-\widehat{\G}\Vert^2\1_{\{\widehat{\G}_a\neq\widehat{\G}\}}\right)\leq a^2\P(\widehat{\G}_a\neq\widehat{\G})
$$
since
\begin{align*}
\Vert\widehat{\G}_a-\widehat{\G}\Vert^2&= \left\Vert\sum_{j=1}^p[\max(\eta_{j,n},a)-\eta_{j,n}]{\bf v}_{jn} {\bf v}_{jn}' \right\Vert^2\nonumber\\
&\leq  \sup_{j=1,...,p}\vert\max(\eta_{j,n},a)-\eta_{j,n}\vert^2\nonumber\\
&\leq  a^2.
\end{align*}
Moreover, since $a < \eta_p = \left\|\G^{-1}\right\|^{-1}$ and by Chebychev inequality, we can bound
\begin{align*}
\P(\widehat{\G}_a\neq\widehat{\G})&= \P(\eta_{p,n}<a)\nonumber\\
&\leq \P \left(\vert\eta_{p,n}-\eta_p\vert\geq\frac{| \eta_p- a|}{2}\right), \\
& \leq  \frac{4}{ (\eta_p- a)^2} \E_p \left(  \vert \eta_{p,n}-\eta_p\vert^2 \right) \\
& \leq    \frac{4}{ (\eta_p- a)^2} \E_p\left(\Vert\widehat{\G}-\G\Vert^2\right),
\end{align*}
because it is known that  the eigenvalue map is Lipschitzian for symmetric matrices   (see \cite{Bhatia96}, Chapter 3). 
This means that for two $p \times p$ symmetric matrices ${\bf A}$ and ${\bf B},$ with eigenvalues $\eta_1({\bf A}) \geq \eta_2({\bf A}) \geq \cdots \geq  \eta_p({\bf A})$ (resp. $\eta_1({\bf B}) \geq \cdots \geq \eta_p({\bf B})$), we have 
$$
\max_{j \in \{ 1, \ldots, p \} } \left| \eta_j({\bf A}) - \eta_j({\bf B}) \right| \leq \left\|{\bf A}-{\bf B}\right\| .
$$
%
%
%

Hence, for some constant $C$
\begin{align}
\E_p\left(\Vert\widehat{\G}_a-\G\Vert^2\right)
&\leq 3\E_p\left(\Vert\widehat{\G}-\G\Vert^2\right)+2a^2 \P(\widehat{\G}_a\neq \widehat{\G}) \nonumber \\
 & \leq  \frac{C}{n}.
\label{equation 4}
\end{align}
Combining   (\ref{equation 1}), (\ref{esp_ga_g}), (\ref{equation 2}) and (\ref{equation 4}), the proof is complete.
\end{Proof}


\begin{Lem}
Under assumptions (A1), (A2) and (A4), there is a constant $C$ such that, for all $n,$  
\begin{align*}
n \ \E_p\left\Vert\frac{1}{N}\sum_{k\in U} \left(\frac{\1_k}{\pi_k}-1 \right)\mathbf{x}_k\right\Vert^2 &\leq C.
\end{align*}
\label{esperancecarrex}
\end{Lem}
\begin{Proof}. Expanding the square norm, we have
\begin{align*}
n\E_p\left\Vert\frac{1}{N}\sum_{k\in U}\alpha_k\mathbf{x}_k\right\Vert^2&= n\E_p\left(\frac{1}{N^2}\sum_{k\in U}\sum_{l\in U}\alpha_{k}\alpha_{l}\mathbf{x}_{k}'\mathbf{x}_{l}\right)\nonumber\\
&\leq  \frac{n}{N^2}\sum_{k\in U}\sum_{l\in U}\left\vert\frac{\Delta_{kl}}{\pi_k\pi_l}\right| \mathbf{x}_{k}'\mathbf{x}_{l} \nonumber\\
&\leq \left[ \frac{n}{N}\frac{1}{\lambda}+\frac{1}{\lambda^2}n\max_{k\neq l\in U}\vert\Delta_{kl}\vert\right]\frac{1}{N}\sum_{k\in U}\Vert\mathbf{x}_k\Vert^2
\end{align*}
and the result follows with hypotheses (A1), (A2) and (A4). 
\end{Proof}

 \begin{Lem}
Under  assumptions  (A2)-(A5), we have 
\begin{itemize}
\item[i)]
$\Vert\tilde{\boldsymbol\beta}(t)-\tilde{\boldsymbol\beta}(r)\Vert^2\leq a^{-2}C_3C_4\vert t-r\vert^{2\beta}.$
\item[ii)] $\Vert\widehat{\boldsymbol\beta}_a(t)-\widehat{\boldsymbol\beta}_a(r)\Vert^2\leq \frac{a^{-2}}{\lambda^2}C_3C_4\vert t-r\vert^{2\beta}.$
\end{itemize}
\label{majoration norme beta}
\end{Lem}

\begin{Proof} 
For i), we just need to remark that, under hypotheses (A3), (A4)  and (A5),
\begin{align}
\Vert\tilde{\boldsymbol\beta}(t)-\tilde{\boldsymbol\beta}(r)\Vert^2&= \left\Vert \G^{-1}\frac{1}{N}\sum_{k\in U}\mathbf{x}_k(Y_k(t)-Y_k(r)\right\Vert^2\nonumber\\
&\leq \Vert \G^{-1}\Vert^2\left(\frac{1}{N}\sum_{k\in U}\Vert\mathbf{x}_k\Vert^2\right)\left(\frac{1}{N}\sum_{k\in U}(Y_k(t)-Y_k(r))^2\right)\nonumber\\
&\leq a^{-2}C_4C_3\vert t-r\vert^{2\beta} . \nonumber
\end{align}

The proof of point ii) is similar, but also requires the use of lower bounds  on the first order inclusion probabilities (assumption (A2)),

\begin{align*}
\Vert\widehat{\boldsymbol\beta}_a(t)-\widehat{\boldsymbol\beta}_a(r)\Vert^2&= \left\Vert \widehat{\G}_a^{-1}\frac{1}{N}\sum_{k\in U}\frac{\1_k}{\pi_k}\mathbf{x}_k(Y_k(t)-Y_k(r)\right\Vert^2\nonumber\\
&\leq \frac{1}{\lambda^2}\Vert \widehat{\G}_a^{-1}\Vert^2\left(\frac{1}{N}\sum_{k\in U}\Vert\mathbf{x}_k\Vert^2\right)\left(\frac{1}{N}\sum_{k\in U}( Y_k(t)-Y_k(r))^2\right)\nonumber\\
&\leq a^{-2}\frac{1}{\lambda^2}C_4C_3\vert t-r\vert^{2\beta} . \nonumber
\end{align*}

\end{Proof}

The following Lemma states the pointwise mean square convergence for any fixed value of $t \in [0,T].$ 
\begin{Lem}
Suppose that assumptions (A1)-(A5) hold. Then, there is a positive constant $\zeta_1$ such that, for all $t \in [0,T],$
\begin{align*}
n\E_p \left(\Vert\widehat{\boldsymbol\beta}_a(t)-\tilde{\boldsymbol\beta}(t)\Vert^2\right) &\leq \zeta_1.
\end{align*}
\label{espbetatinstantt}
\end{Lem}

\begin{Proof}. The demonstration is similar to the proof of Lemma  \ref{esperance diff beta}
and is thus omitted.
\end{Proof}

\begin{Lem}
Suppose that assumptions (A1)-(A5) hold. Then, there is a positive constant  $\zeta_2$ such that
\begin{align*}
n\E_p\left( \Vert\widehat{\boldsymbol\beta}_a(t)-\tilde{\boldsymbol\beta}(t)-\widehat{\boldsymbol\beta}_a(r)+\tilde{\boldsymbol\beta}(r)\Vert^2\right) &\leq\zeta_2\vert t-r\vert^{2\beta}.
\end{align*}
\label{esperance diff beta}
\end{Lem}

\begin{Proof}. 
A direct decomposition leads to
\begin{align}
 &  n\Vert\widehat{\boldsymbol\beta}_a(t)-\tilde{\boldsymbol\beta}(t)-\widehat{\boldsymbol\beta}_a(r)+\tilde{\boldsymbol\beta}(r)\Vert^2\nonumber\\
&\leq  \Biggl\Vert (\widehat{\G}_a^{-1}-\G^{-1})\frac{1}{N}\sum_{k\in U}\frac{\1_k}{\pi_k}\mathbf{x}_k(Y_k(t)-Y_k(r)) + \G^{-1}\frac{1}{N}\sum_{k\in U}\left(\frac{\1_k}{\pi_k}-1\right)\mathbf{x}_k(Y_k(t)-Y_k(r))\Biggr\Vert^2 \nonumber\\
&\leq 2A_{1N}^2+2A_{2N}^2,
\label{diff beta_hat beta}
\end{align}
where $A_{1N}^2=n\Vert \widehat{\G}_a^{-1}-\G^{-1}\Vert^2\Bigl\Vert\frac{1}{N}\sum_{k\in U}\frac{\1_k}{\pi_k}\mathbf{x}_k(Y_k(t)-Y_k(r))\Bigr\Vert^2 $
 and \\  $A_{2N}^2=n \Vert \G^{-1}\Vert^2\Bigl\Vert\frac{1}{N}\sum_{k\in U}\alpha_k\mathbf{x}_k(Y_k(t)-Y_k(r))\Bigr\Vert^2.$
Using now assumptions (A2)-(A4) and the Cauchy-Schwarz inequality, we get
\begin{align}
A_{1N}^2 &\leq  n\Vert\widehat{\G}_a^{-1}-\G^{-1}\Vert ^2\left(\frac{1}{\lambda^2}\frac{1}{N}\sum_{k\in U}\Vert\mathbf{x}_k\Vert ^2\right)\left(\frac{1}{N}\sum_{k\in U} ( Y_k(t)-Y_k(r)) ^2 \right)\nonumber\\
&\leq n\Vert\widehat{\G}_a^{-1}-\G^{-1}\Vert ^2\frac{1}{\lambda^2}C_3C_4\vert t-r\vert^{2\beta}. \nonumber
\end{align}
Using now Lemma \ref{norme Ga-G}, we can bound
\begin{align}
\E_p(A_{1N}^2)&\leq C \vert t- r\vert^{2\beta},
\label{esperance A1N}
\end{align}
for some constant $C.$ Now,  with assumptions (A1)-(A5) and following the same arguments as in the proof of Lemma \ref{esperancecarrex}, we also have
\begin{align}
\E_p(A_{2N}^2)&\leq n\Vert \G^{-1}\Vert^2 \E_p \left( \left\| \frac{1}{N}\sum_{k\in U}\alpha_k\mathbf{x}_k( Y_k(t)-Y_k(r))\right\|^2 \right) \nonumber\\
&\leq  \left(\frac{n}{N}\frac{1}{\lambda}+\frac{n\max_{k\neq l\in U}\vert\Delta_{kl}\vert}{\lambda^2}\right)C_3C_4a^{-2}\vert t-r\vert^{2\beta} \leq C\vert t-r\vert^{2\beta}.
\label{esperance A2N}
\end{align}
for some positive constant $C.$
%
%
%
Combining (\ref{diff beta_hat beta}), (\ref{esperance A1N}) and  (\ref{esperance A2N}), the result is proved. 

\end{Proof}

\subsection{Proof of Proposition \ref{convergencebeta} and Proposition \ref{convergence mu}}

The proof of Proposition \ref{convergencebeta} is omitted. It is analogous to the proof of Proposition \ref{convergence mu}, which is given below.
The different steps are similar to the proof of Proposition 1 in  \cite{cardot_josserand_2011}.

\medskip
Let us decompose, for $t \in [0,T],$
\begin{equation}
\displaystyle\sup_{t\in[0,T]}\vert\widehat{\mu}_{\textrm{MA},d}(t)-\mu(t)\vert\leq \displaystyle\sup_{t\in[0,T]}\vert\widehat{\mu}_{\textrm{MA},d}(t)-\widehat{\mu}_{\textrm{MA},a}(t)\vert+\sup_{t\in[0,T]}\vert\widehat{\mu}_{\textrm{MA},a}(t)-\mu(t)\vert
\label{erreur approximation + estimation}
\end{equation}
and study each term at the right-hand side of the inequality separately.

\subsubsection*{Step 1. The interpolation error $\sup_{t\in[0,T]}\vert\widehat{\mu}_{\textrm{MA},d}(t)-\widehat{\mu}_{\textrm{MA},a}(t)\vert.$}

Consider $t\in[t_i,t_{i+1})$ and  write 
\begin{align*}
\vert\widehat{\mu}_{\textrm{MA},d}(t)-\widehat{\mu}_{\textrm{MA},a}(t)\vert&\leq \vert\widehat{\mu}_{\textrm{MA},a}(t_i)-\widehat{\mu}_{\textrm{MA},a}(t)\vert+\vert\widehat{\mu}_{\textrm{MA},a}(t_{i+1})-\widehat{\mu}_{\textrm{MA},a}(t_i)\vert.
\end{align*}
\noindent Under assumptions (A2)-(A5) and using Lemma \ref{majoration norme beta}, ii), we get 
\begin{align*}
\vert\widehat{\mu}_{\textrm{MA},a}(t)-\widehat{\mu}_{\textrm{MA},a}(r)\vert&\leq \left\vert\frac{1}{N}\sum_{k\in U}\alpha_k\mathbf{x}_k'(\widehat{\boldsymbol{\beta}}_a(t)-\widehat{\boldsymbol{\beta}}_a(r))\right\vert+\frac{1}{N}\sum_{k\in s}\frac{\vert Y_k(t)-Y_k(r)\vert}{\pi_k}\nonumber\\
&\leq \left(1+\frac{1}{\lambda}\right)\sqrt{C_4}\Vert\widehat{\boldsymbol{\beta}}_a(t)-\widehat{\boldsymbol{\beta}}_a(r)\Vert+\frac{1}{\lambda}\left(\frac{1}{N}\sum_{k\in U}(Y_k(t)-Y_k(r))^2\right)^{1/2}\nonumber\\
&\leq \left((1+\lambda^{-1})C_4a^{-1}+1\right)\lambda^{-1}\sqrt{C_3}\vert t-r\vert^\beta.\nonumber
\end{align*}
So, there is a positive constant $C$ such that 
\begin{align*}
\vert\widehat{\mu}_{\textrm{MA},a}(t)-\widehat{\mu}_{\textrm{MA},a}(r)\vert &\leq C\vert t-r\vert^\beta 
\end{align*}
and consequently, 
\begin{align*}
\vert\widehat{\mu}_{\textrm{MA},d}(t)-\widehat{\mu}_{\textrm{MA},a}(t)\vert&\leq C[\vert t_i-t\vert^\beta+\vert t_{i+1}-t_i\vert^\beta]\nonumber\\
&\leq 2C\vert t_{i+1}-t_i\vert^\beta.
\end{align*}
\noindent Hence, since by hypothesis, $\lim_{N\rightarrow \infty}\max_{i=\{1,...,d_{N-1}\}}\vert t_{i+1}-t_i\vert^{\beta}=o(n^{-1/2})$, we have

\begin{align}
\sup_{t \in [0,T]} \ \sqrt{n}\vert\widehat{\mu}_{\textrm{MA},d}(t)-\widehat{\mu}_{\textrm{MA},a}(t)\vert &= o(1).
\label{convergencemuinterpolee}
\end{align}

\subsubsection*{Step 2. The estimation error $\sup_{t \in [0,T]} \vert\widehat{\mu}_{{MA},a}(t)-\mu(t)\vert.$}

We use the following decomposition:
\begin{align}
\sup_{t\in [0,T]}\left| \widehat{\mu}_{\textrm{MA},a}(t)-\mu(t)\right| & \leq   \vert\widehat{\mu}_{\textrm{MA},a}(0)-\mu(0)\vert+\mbox{sup}_{r,t\in [0,T]}\vert\widehat{\mu}_{\textrm{MA},a}(t)-\mu(t)-\widehat{\mu}_{\textrm{MA},a}(r)+  \mu(r)\vert .\nonumber\\
\label{step2.1}
\end{align}
Writing,
 \begin{align*}
 \widehat{\mu}_{\textrm{MA},a}(0) - \mu(0) & =    \frac{1}{N} \sum_{k\in U}\alpha_kY_k(0)-\frac{1}{N} \sum_{k\in U}\alpha_k\widehat{Y}_k(0)\\
  & =  \frac{1}{N} \sum_{k\in U} \alpha_kY_k(0) -\frac{1}{N^2} \sum_{k\in U} \alpha_k{\bf x}_k'\widehat{\bf G}_a^{-1}\sum_{l\in s}\frac{{\bf x}_l Y_l(0)}{\pi_l}
 \end{align*}
we directly get, with hypotheses A1-A3  and with similar arguments as in the proof of Lemma   \ref{esperancecarrex},
that  for some constant $C,$
\begin{align}
\E_p \left( \widehat{\mu}_{\textrm{MA},a}(0) - \mu(0) \right)^2 & \leq \frac{C}{n}.
\label{step1.2}
\end{align} 

\noindent  The second term at the right-hand side  in (\ref{step2.1}) is dealt with using maximal inequalities. More exactly, we use Corollary 2.2.5 in \cite{vdVWellener2000}.  Consider for this, the  Orlicz norm of some random variable $X$ which is defined as follows
\begin{align*}
||X||_{\psi} &=\sqrt{\E_p(\psi(X))}.
\end{align*}
For the particular case $\psi(u)=u^2,$  the Orlicz norm  is simply  the well-known $L^2$ norm, $||X||_{\psi}=\sqrt{\E_p(X^2)}.$ Let us introduce for $(r,t)  \in [0,T]^2,$ the semimetric $d(r,t)$ defined by 
\begin{align*}
d^2(r,t) &=  \left\| \sqrt{n}\vert\widehat{ \mu}_{\textrm{MA},a}(t)-\mu(t)-\widehat{\mu}_{\textrm{MA},a}(r)-\mu(r)\vert \right\|^2_{\psi}\\
 & =  n\E_p \left(\vert\widehat{\mu}_{\textrm{MA},a}(t)-\mu(t)-\widehat{\mu}_{\textrm{MA},a}(r)+\mu(r)\vert^2 \right)
\end{align*}
and consider  $D(\epsilon, d),$ the packing number, which is defined as  the maximum number of points in $[0,T]$ whose distance between each pair is strictly larger than $\epsilon.$ Then,
Corollary 2.2.5 in \cite{vdVWellener2000} states that   there is a constant $K>0$ such that 
\begin{align}
\left|\left|\sup_{(r,t) \in [0,T]^2 }\sqrt{n}\vert\widehat{ \mu}_{\textrm{MA},a}(t)-\mu(t)-\widehat{\mu}_{\textrm{MA},a}(r)-\mu(r)\vert \right|\right|_{\psi} 
 &\leq  K \int_0^{T} \psi^{-1}( D(\epsilon,d))d\epsilon.
 \label{maxinequalitymu}
\end{align}
We show below that there is a constant $C$ such that $
d^2(r,t)\leq C\vert t-r\vert^{2\beta}$ and thus, since $\beta>1/2,$ the integral at the right-hand side of (\ref{maxinequalitymu}) is finite. 

Let us first decompose
\begin{align}
d^2(r,t)&\leq 2d_1^2(r,t)+2d_2^2(r,t)
\label{distance mu hat - mu}
\end{align}
where 
\begin{align*}
d_1^2(r,t) &=n\E_p(\vert\widehat{\mu}_{\textrm{MA},a}(t)-\tilde{\mu}(t)-\widehat{\mu}_{\textrm{MA},a}(r)+\tilde{\mu}(r)\vert^2)
\end{align*}
and 
\begin{align*}
d_2^2(r,t) &=n\E_p(\vert\tilde{\mu}(t)-\mu(t)-\tilde{\mu}(r)+\mu(r)\vert^2).
\end{align*}



By assumptions (A2)-(A4) and Lemma \ref{esperance diff beta}, we can bound, for some constant $C,$
\begin{align}
d_1^2(r,t)
&\leq  \E_p\left\lbrace n\left\Vert\frac{1}{N}\sum_{k\in U}\alpha_k  \mathbf{x}_k\right\Vert^2\Vert\widehat{\boldsymbol\beta}_a(t)-\tilde{\boldsymbol\beta}(t)-\widehat{\boldsymbol\beta}_a(r)+\tilde{\boldsymbol\beta}(r)\Vert^2\right\rbrace\nonumber\\
& \leq   \left(1+\frac{1}{\lambda}\right)^2C_4\zeta_2\vert t-r\vert^{2\beta}:=C \vert t-r\vert^{2\beta}.
\label{erreur estimation}
\end{align}
 Considering now $d_2(r,t),$ we have
\begin{align}
d_2^2(r,t)&= n\E_p\left[ \frac{1}{N}\sum_{k\in U} \alpha_k\left[Y_k(t)-Y_k(r)-\mathbf{x}_k'(\tilde{\boldsymbol\beta}(t)-\tilde{\boldsymbol\beta}(r))\right]\right]^2\nonumber\\
& \leq  2\E_p(A_N^2)+2\E_p(B_N^2)
\end{align}
where $A_N^2=n\left( \frac{1}{N}\sum_{k\in U} \alpha_k\left[ Y_k(t)-Y_k(r)\right]\right)^2$ and $B_N^2=n\left( \frac{1}{N}\sum_{k\in U} \alpha_k\mathbf{x}_k'(\tilde{\boldsymbol\beta}(t)-\tilde{\boldsymbol\beta}(r))\right)^2.$
With hypotheses (A1)-(A3), one can easily obtain that there is a positive constant $C$ such that
\begin{align}
\E_p(A_N^2)
&\leq  C \vert t-r\vert^{2\beta}
\label{equation BN}
\end{align} 
and thanks to Lemma \ref{esperancecarrex} and to Lemma  \ref{majoration norme beta}, we can bound
\begin{align}
\E_p(B_N^2)
&\leq \E_p\left[n\left\Vert \frac{1}{N}\sum_{k\in U} \alpha_k\mathbf{x}_k\right\Vert^2\right]\Vert\tilde{\boldsymbol\beta}(t)-\tilde{\boldsymbol\beta}(r)\Vert^2\nonumber\\
&\leq C \vert t-r\vert^{2\beta}.
\label{equation CN}
\end{align}
Combining now  (\ref{equation BN}) and (\ref{equation CN})
 with (\ref{distance mu hat - mu}) and (\ref{erreur estimation}), we get  that 
\begin{align}
d^2(r,t)&\leq C\vert t-r\vert^{2\beta},
\label{majoration d2}
\end{align}
for some constant $C.$

Using now (\ref{majoration d2}), it is clear that the packing number is  bounded as follows:  $D(\epsilon, d) = O(\epsilon^{-1/\beta}).$ Consequently, the integral at the right-hand side of (\ref{maxinequalitymu}) is finite when $\beta>1/2.$ Inserting (\ref{step1.2}) and (\ref{maxinequalitymu}) in (\ref{step2.1}), the proof of step 2 is complete.

\subsection{Proof of the consistency of the covariance function}


We first prove that for any $(r,t) \in [0,T]^2,$ the estimator $\widehat{\gamma}_{\textrm{MA},d}(r,t)$ of the covariance function converges to $\gamma_{\textrm{MA}}(r,t).$ 

Then we prove  the uniform convergence of the variance estimator $\widehat{\gamma}_{\textrm{MA},d}(t,t)$ by showing its convergence in distribution to zero in the space of continuous functions.  The proof is decomposed into two classical steps (see for example Theorem 8.1 in \cite{Billingsley1968}). We first show the pointwise convergence, by considering the convergence of all finite linear combinations, and then we check that the sequence is tight by bounding the increments.

\subsubsection*{Step 1. Pointwise convergence}

We want to show, that for each $(t,r)\in[0,T]^2$, we have when $N$ tends to infinity,
\begin{align*}
n \E_p\left\lbrace \mid\widehat{\gamma}_{\textrm{MA},d}(r,t)-\gamma_{\textrm{MA}}(r,t)\mid\right\rbrace &\rightarrow  0.
\end{align*}
Let us decompose
\begin{align*}
n(\widehat{\gamma}_{\textrm{MA},d}(r,t)- {\gamma}_{\textrm{MA}}(r,t))&= n(\widehat{\gamma}_{\textrm{MA},d}(r,t)- \widehat{\gamma}_{\textrm{MA},a}(r,t))+n(\widehat{\gamma}_{\textrm{MA},a}(r,t)- \gamma_{\textrm{MA}}(r,t))
\end{align*}
where $\widehat{\gamma}_{\textrm{MA},a}(r,t)$ is defined by 
\begin{align*}
\widehat{\gamma}_{\textrm{MA},a}(r,t) &=\frac{1}{N^2}\sum_{k,l \in s}\frac{\Delta_{kl}}{ \pi_{kl}}\frac{Y_k(r)-\widehat{Y}_{k,a}(r)}{\pi_k}\cdot\frac{Y_l(t)-\widehat{Y}_{l,a}(t)}{\pi_l}
\end{align*}
We study separately the interpolation and the estimation errors.

\subsubsection*{Interpolation error}

Let us suppose that $t\in[t_i,t_{i+1})$, $r \in[t_{i'},t_{i'+1}).$ We have $n(\widehat{\gamma}_{\textrm{MA},d}(r,t)- \widehat{\gamma}_{\textrm{MA},a}(r,t))\leq A+B,$ with 
\begin{align*}
A & =  \frac{n}{N^2} \sum_{k,l \in U}\frac{\vert\Delta_{kl}\vert}{ \pi_{kl}\pi_k\pi_l}|(Y_{k,d}(r)-Y_k(r))(Y_{l,d}(t)-Y_l(t))\\
&+(Y_{k,d}(r)-Y_k(r))(Y_l(t)-\widehat Y_{l,d}(t))+(Y_k(r)-\widehat Y_{k,d}(r))(Y_{l,d}(t)-Y_l(t))|
\end{align*}
and
\begin{align*} 
B&=\frac{n}{N^2} \sum_{k,l \in U}\frac{\vert\Delta_{kl}\vert}{ \pi_{kl}\pi_k\pi_l}\Biggl\vert\left( Y_k(r)-\widehat{Y}_{k,d}(r) \right)\left( Y_l(t)-\widehat{Y}_{l,d}(t) \right)-\left( Y_k(r)-\widehat{Y}_{k,a}(r) \right)\left( Y_l(t)-\widehat{Y}_{l,a}(t) \right)\Biggr\vert\\
 & =  \frac{n}{N^2} \sum_{k,l \in U}\frac{\vert\Delta_{kl}\vert}{ \pi_{kl}\pi_k\pi_l}\Biggl\vert Y_k(r)(\widehat Y_{l,a}(t)-\widehat Y_{l,d}(t))+Y_l(t)(\widehat Y_{k,a}(r)-\widehat Y_{k,d}(r))+\widehat Y_{k,d}(r)\widehat Y_{l,d}(t)-\widehat Y_{k,a}(r)\widehat Y_{l,a}(t)\Biggr\vert .
\end{align*}
For $t\in [t_i,t_{i+1}],$ we can write
$$
|Y_{l,d}(t)-Y_l(t)|\leq |Y_l(t_i)-Y_l(t)|+|Y_l(t_{i+1})-Y_l(t_i)|
$$
and
$$
|\widehat Y_{l,a}(t)-\widehat Y_{l,d}(t)|\leq |\widehat Y_{l,a}(t)-\widehat Y_{l,a}(t_i)|+|\widehat Y_{l,a}(t_{i+1})-\widehat Y_{l,d}(t_i)|
$$
We have that $\frac{1}{N}\sum_{l\in U}(Y_{l,d}(t)- Y_{l}(t))^2\leq C[|t_i-t|^{2\beta}+|t_{i+1}-t_i|^{2\beta}]$ and $\frac{1}{N}\sum_{l\in U}(Y_l(t)-\widehat{Y}_{l,d}(t))^2=O(1). $
Thanks to Lemma \ref{majoration norme beta}, we can bound 
\begin{eqnarray}
\vert\widehat{Y}_{l,a}(t_{i})-\widehat{Y}_{l,a}(t)\vert&\leq&C_4a^{-1}\frac{1}{\lambda}C_3^{1/2}\vert t_{i}-t\vert^\beta\leq C_4a^{-1}\frac{1}{\lambda}C_3^{1/2}\vert t_{i+1}-t_i\vert^\beta. \nonumber
\end{eqnarray}
Under the assumption on the grid of discretization points, one can get after some algebra that
\begin{eqnarray*}
n\vert\widehat{\gamma}_{\textrm{MA},d}(r,t)- \widehat{\gamma}_{\textrm{MA},a}(r,t)\vert=o(1).
\end{eqnarray*}

\subsubsection*{Estimation error}
Consider now, 
\begin{align}
n(\widehat{\gamma}_{\textrm{MA},a}(r,t)- \gamma_{\textrm{MA}}(r,t))&=\frac{n}{N^2}\sum_U\sum_U\frac{\Delta_{kl}}{\pi_k\pi_l}\left(\frac{\1_{kl}}{\pi_{kl}}-1\right)[Y_k(t)-\tilde{Y}_k(t)][Y_l(r)-\tilde{Y}_l(r)]\nonumber\\
&+ \frac{n}{N^2}\sum_{k\in U}\sum_{l\in U}\frac{\Delta_{kl}}{\pi_k\pi_l}\frac{\1_{kl}}{\pi_{kl}}[Y_k(t)-\tilde{Y}_k(t)][\tilde{Y}_l(r)-\widehat{Y}_{l,a}(r)]\nonumber\\
&+ \frac{n}{N^2}\sum_{k\in U}\sum_{l\in U}\frac{\Delta_{kl}}{\pi_k\pi_l}\frac{\1_{kl}}{\pi_{kl}}[\tilde{Y}_k(t)-\widehat{Y}_{k,a}(t)][Y_l(r)-\tilde{Y}_l(r)]\nonumber\\
&+ \frac{n}{N^2}\sum_{k\in U}\sum_{l\in U}\frac{\Delta_{kl}}{\pi_k\pi_l}\frac{\1_{kl}}{\pi_{kl}}[\tilde{Y}_k(t)-\widehat{Y}_{k,a}(t)][\tilde{Y}_l(r)-\widehat{Y}_{l,a}(r)]\nonumber\\
&:=A_1(r,t)+A_2(r,t)+A_3(r,t)+A_4(r,t).
\label{def:decompGma}
\end{align}
Let us define  $\tilde{e}_k(t)=Y_k(t)-\tilde{Y}_k(t)$ and first show that $\E_p(A_1(r,t)^2)\rightarrow0$ when $N\rightarrow\infty$.

\begin{align}
\E_p(A_1(r,t)^2)&=\E_p\left[\frac{n^2}{N^4}\sum_{k,l\in U}\sum_{k',l'\in U}\frac{\Delta_{kl}}{\pi_k\pi_l}\left(\frac{\1_{kl}}{\pi_{kl}}-1\right)\frac{\Delta_{k'l'}}{\pi_{k'}\pi_{l'}}\left(\frac{\1_{k'l'}}{\pi_{k'l'}}-1\right)\tilde{e}_k(t)\tilde{e}_l(r)\tilde{e}_{k'}(t)\tilde{e}_{l'}(r)\right]\nonumber\\
&=\E_p\left[\frac{n^2}{N^4}\sum_{k\in U}\sum_{k'\in U}\frac{1-\pi_k}{\pi_k}\left(\frac{\1_{k}}{\pi_{k}}-1\right)\frac{1-\pi_{k'}}{\pi_{k'}}\left(\frac{\1_{k'}}{\pi_{k'}}-1\right)\tilde{e}_k(t)\tilde{e}_k(r)\tilde{e}_{k'}(t)\tilde{e}_{k'}(r)\right]\nonumber\\
&+2\E_p\left[\frac{n^2}{N^4}\sum_{k\in U}\sum_{k'\neq l'\in U}\frac{1-\pi_k}{\pi_k}\left(\frac{\1_{k}}{\pi_{k}}-1\right)\frac{\Delta_{k'l'}}{\pi_{k'}\pi_{l'}}\left(\frac{\1_{k'l'}}{\pi_{k'l'}}-1\right)\tilde{e}_k(t)\tilde{e}_k(r)\tilde{e}_{k'}(t)\tilde{e}_{l'}(r)\right]\nonumber\\
&+\E_p\left[\frac{n^2}{N^4}\sum_{k\neq l\in U}\sum_{k'\neq l'\in U}\frac{\Delta_{kl}}{\pi_k\pi_l}\left(\frac{\1_{kl}}{\pi_{kl}}-1\right)\frac{\Delta_{k'l'}}{\pi_{k'}\pi_{l'}}\left(\frac{\1_{k'l'}}{\pi_{k'l'}}-1\right)\tilde{e}_k(t)\tilde{e}_l(r)\tilde{e}_{k'}(t)\tilde{e}_{l'}(r)\right]\nonumber\\
&:=a_1+a_2+a_3.
\label{d1}
\end{align}
The hypotheses on the moments of the inclusion probabilities and Lemma  \ref{somme e tilde puissance 2} give us
\begin{align*}
a_1
&\leq\left(\frac{n^2}{N^3}\frac{1}{\lambda^3}+\frac{n^2}{N^2}\frac{\max_{k\neq k'\in U}\vert\Delta_{kk'}\vert}{\lambda^4}\right)\zeta_4
\end{align*}
as well as 
\begin{align*}
a_3
&\leq\frac{C}{N}+ \frac{(n\max_{k\neq l\in U}\vert\Delta_{kl}\vert)^2}{\lambda^4\lambda*^2}\max_{(k,l,k',l')\in D_{4,n}}\vert \E_p\{(\1_{kl}-\pi_{kl})(\1_{k'l'}-\pi_{k'l'})\}\vert\zeta_5
\end{align*}
so that  $a_1\rightarrow0$ and  $a_3\rightarrow0$ when $N\rightarrow\infty$.
Then, the Cauchy-Schwarz inequality allows us to get that 
  $a_2\rightarrow0$ when $N\rightarrow\infty$ and 
 $\E_p(A_1(r,t)^2)\rightarrow0$ when $N\rightarrow\infty$.

Let us show now that  $\E_p(\vert A_4(r,t)\vert)\rightarrow0 $ when  $N\rightarrow\infty$.
With Lemma \ref{espbetatinstantt}, and assumptions (A1)-(A5), we have
\begin{align*}
\E_p(\vert A_4(r,t)\vert)
&\leq n \E_p\Biggl(\frac{1}{N^2}\sum_{k\in U}\sum_{l\in U}\frac{\vert\Delta_{kl}\vert}{\pi_k\pi_l}\frac{1}{\pi_{kl}}\Vert\mathbf{x}_k\Vert\Vert\mathbf{x}_l\Vert\Vert \tilde{\boldsymbol\beta}(t)-\widehat{\boldsymbol\beta}_a(t)\Vert\Vert\tilde{\boldsymbol\beta}(r)-\widehat{\boldsymbol\beta}_a(r)\Vert\Biggr)\nonumber\\
&\leq \frac{1}{n}\left[\frac{n}{\lambda^2N}+\frac{n\max_{k\neq l\in U}\vert\Delta_{kl}\vert}{\lambda^2\lambda*}\right]C_4\zeta_1
\end{align*}
so that $\E_p(\vert A_4(r,t)\vert)\rightarrow0$ when $N\rightarrow\infty$.

In a similar way, we can bound 
 $\E_p(\vert A_2(r,t)\vert)$ as follows,
\begin{align*}
\E_p(\vert A_2(r,t)\vert)&\leq\frac{n}{N^2}\sum_{k\in U}\sum_{l\in U}\frac{\vert\Delta_{kl}\vert}{\pi_k\pi_l}\frac{1}{\pi_{kl}}\E_p\vert\tilde{e}_k(t)\widehat{\tilde{e}}_l(r)\vert\nonumber\\
&\leq \frac{n}{N^2}\sum_{k\in U}\sum_{l\in U}\frac{\vert\Delta_{kl}\vert}{\pi_k\pi_l}\frac{\Vert\mathbf{x}_l\Vert}{\pi_{kl}}\vert Y_k(t)-\tilde{Y}_k(t) \vert\cdot \E_p(\Vert\tilde{\boldsymbol\beta}(r)-\widehat{\boldsymbol\beta}_a(r)\Vert)\nonumber\\
&\leq \left(\frac{\sqrt{n}}{\lambda^2N}+\frac{\sqrt{n}\max_{k\neq l\in U}\vert\Delta_{kl}\vert}{\lambda^2\lambda^*}\right)C_4^{1/2}\zeta_1^{1/2}\frac{1}{N}\sum_{k\in U}\vert Y_k(t)-\tilde{Y}_k(t)\vert,
\end{align*}
where $\widehat{\tilde{e}}_k(t)=\tilde{Y}_k(t)-\widehat{Y}_{k,a}(t)=\mathbf{x}'_k(\tilde{\boldsymbol{\beta}}(t)-\hat{\boldsymbol{\beta}}_a(t)).$
Thus, there is constant $C$ such that,
\begin{align*}
\E_p(\vert A_2(r,t)\vert)
&\leq \frac{C}{\sqrt{n}}
\end{align*}
and $\E_p(\vert A_2(r,t)\vert)\rightarrow 0$ when $N\rightarrow\infty.$ We can show in a similar way that 
$\E_p(\vert A_3(r,t)\vert)\rightarrow 0$ when $N\rightarrow\infty$.

Finally, we have that for all $(r,t) \in [0,T]^2,$
\begin{align}
n \E_p\left\lbrace \mid\widehat{\gamma}_{\textrm{MA},a}(r,t)-\gamma_{\textrm{MA}}(r,t)\mid\right\rbrace &\rightarrow  0, \quad \text{ when } N\rightarrow\infty.
\label{Convergence ponctuelle cov}
\end{align}






\subsubsection*{Step 2. Uniform convergence of the variance estimator}

The  pointwise convergence of the variance function proved in the previous step clearly implies the convergence of all finite linear combinations : for all $p \in \{1,2, \ldots \}$, for all $(c_1, \ldots, c_p) \in \mathbb{R}^p$ and for all $(t_1, \ldots, t_p) \in [0,T]^p,$ we have 
\begin{align}
\sum_{\ell =1}^p c_\ell \ n \left( \widehat{\gamma}_{\textrm{MA},a}(t_\ell,t_\ell)  - \gamma_{\textrm{MA}}(t_\ell,t_\ell) \right)& \rightarrow   0
\end{align}
in probability as $N$ tends to infinity. Thus, we deduce with the Cramer-Wold device that the vector $n \left( \widehat{\gamma}_{\textrm{MA},a}(t_1,t_1)  - \gamma_{\textrm{MA}}(t_1,t_1), \ldots, \widehat{\gamma}_{\textrm{MA},a}(t_p,t_p)  - \gamma_{\textrm{MA}}(t_p,t_p) \right)$ converges in distribution to $0$ (in $\mathbb{R}^p$).

We need now to prove that the sequence of random functions $\widehat{\gamma}_{\textrm{MA},a}(t,t)$ is tight in $C[0,T]$ by using a bound on its increments.
Let us introduce the following criterion, 
\begin{align*}
d^2_{\gamma}(t,r)&=n^2 \E_p(\vert \widehat{\gamma}_{\textrm{MA},a}(t,t)-\gamma_{\textrm{MA}}(t,t)-\widehat{\gamma}_{\textrm{MA},a}(r,r)+\gamma_{\textrm{MA}}(r,r)\vert^2).
\end{align*}
To conclude we show in the following that $d^2_{\gamma}(t,r) \leq C \vert t-r\vert^{2\beta}$ for a constant $C$ and all $(r,t) \in [0,T]^2.$
Using (\ref{def:decompGma}), the distance is decomposed into four parts.

Let us define $\phi_{kl}(t,r)=\tilde{e}_{k}(t)\tilde{e}_l(t)-\tilde{e}_{k}(r)\tilde{e}_l(r)$ and first consider 
 $d^2_{A_1}=\E_p(\vert A_1(t,t)-A_1(r,r)\vert^2).$
We have
\begin{align}
d^2_{A_1}&= \E_p\left[\frac{n^2}{N^4}\sum_{k\in U}\sum_{k'\in U}\frac{1-\pi_k}{\pi_k}\left(\frac{\1_{k}}{\pi_{k}}-1\right)\frac{1-\pi_{k'}}{\pi_{k'}}\left(\frac{\1_{k'}}{\pi_{k'}}-1\right)\phi_{kk}(t,r)\phi_{k'k'}(t,r)\right]\nonumber\\
& +2\E_p\left[\frac{n^2}{N^4}\sum_{k\in U}\sum_{k'\neq l'\in U}\frac{1-\pi_k}{\pi_k}\left(\frac{\1_{k}}{\pi_{k}}-1\right)\frac{\Delta_{k'l'}}{\pi_{k'}\pi_{l'}}\left(\frac{\1_{k'l'}}{\pi_{k'l'}}-1\right)\phi_{kk}(t,r)\phi_{k'l'}(t,r)\right]\nonumber\\
& +\E_p\left[\frac{n^2}{N^4}\sum_{k\neq l\in U}\sum_{k'\neq l'\in U}\frac{\Delta_{kl}}{\pi_k\pi_l}\left(\frac{\1_{kl}}{\pi_{kl}}-1\right)\frac{\Delta_{k'l'}}{\pi_{k'}\pi_{l'}}\left(\frac{\1_{k'l'}}{\pi_{k'l'}}-1\right)\phi_{kl}(t,r)\phi_{k'l'}(t,r)\right]\nonumber\\
&:= b_1+b_2+b_3
\label{d A1 puissance 2}
\end{align} 
Thanks to Lemma  \ref{somme phi puissance 2}, we get
\begin{align}
b_1&\leq \left(\frac{n^2}{N^3}\frac{1}{\lambda^3}+\frac{n^2}{N^2}\frac{\max_{k\neq k'\in U}\vert\Delta_{kk'}\vert}{\lambda^4}\right)\frac{1}{N}\sum_{k\in U}\vert\phi_{kk}(t,r)\vert^2 \nonumber\\
 &\leq C\vert t-r\vert^{2\beta}
\label{b1}
\end{align}
and
\begin{align}
b_3&\leq
\frac{C}{N}|t-r|^{2 \beta} + \frac{(n\max_{k\neq l\in U}\vert\Delta_{kl}\vert)^2}{\lambda^4\lambda*^2}\max_{(k,l,k',l')\in D_{4,n}}\vert \E_p\{(\1_{kl}-\pi_{kl})(\1_{k'l'}-\pi_{k'l'})\}\vert\left(\frac{1}{N^2}\sum_{k,l\in U}\vert\phi_{kl}(t,r)\vert\right)^2\nonumber\\
&\leq C\vert t-r\vert^{2\beta}.
\label{b3}
\end{align}

The Cauchy-Schwarz inequality together with bounds (\ref{b1}) and (\ref{b3}) allows us to get
$
b_2
\leq C\vert t-r\vert^{2\beta}$
so that
\begin{align}
d^2_{A_1}&\leq  C\vert t-r\vert^{2\beta}.
\label{d2 A1}
\end{align}

Let us bound now 
$d^2_{A_2}=\E_p(\vert A_2(t,t)-A_2(r,r)\vert^2)$ and define 
 $\tilde{\phi}_{kl}(t,r)=\tilde{e}_k(t)\widehat{\tilde{e}}_{l}(t)-\tilde{e}_k(r)\widehat{\tilde{e}}_{l}(r).$
 Thanks to Lemma \ref{somme phi tilde puissance 2}, we get
\begin{align}
d^2_{A_2}&\leq\frac{2n^2}{N^2\lambda^4}\E_p\left(\frac{1}{N}\sum_{k\in U}\tilde{\phi}_{kk}(t,r)\right)^2+\frac{2n^2\max_{k\neq l\in U}\vert\Delta_{kl}\vert^2}{\lambda^4\lambda*^2}\E_p\left(\frac{1}{N^2}\sum_{k,l\in U}\vert\tilde{\phi}_{k,l}(t,r)\vert\right)^2\nonumber\\
&\leq C\vert t-r\vert^{2\beta}.
\label{d2 A2}
\end{align}

Let us study now the last term, $d^2_{A_4}=\E_p(\vert A_4(t,t)-A_4(r,r)\vert^2)$ 
and define $\widehat{\tilde{\phi}}_{kl}(t,r)=\widehat{\tilde{e}}_k(t)\widehat{\tilde{e}}_{l}(t)-\widehat{\tilde{e}}_k(r)\widehat{\tilde{e}}_{l}(r).$
Thanks to  Lemma \ref{somme phi tilde hat puissance 2}, we have
\begin{align}
d^2_{A_4}
&\leq \frac{2n^2}{N^2\lambda^4}\E_p\left(\frac{1}{N}\sum_{k\in U}\widehat{\tilde{\phi}}_{kk}(t,r)\right)^2+\frac{2n^2\max_{k\neq l\in U}\vert\Delta_{kl}\vert^2}{\lambda^4\lambda*^2}\E_p\left(\frac{1}{N^2}\sum_{k,l\in U}\vert\widehat{\tilde{\phi}}_{k,l}(t,r)\vert\right)^2\nonumber\\
&\leq C\vert t-r\vert^{2\beta}.
\label{d2 A4}
\end{align}

Finally, we can deduce, with inequalities (\ref{def:decompGma}), (\ref{d2 A1}), (\ref{d2 A2}) and (\ref{d2 A4}), that
\begin{align}
d^2_{\gamma}(t,r)&= n^2 \E_p(\vert \widehat{\gamma}_{\textrm{MA},a}(t,t)-\gamma_{\textrm{MA}}(t,t)-\widehat{\gamma}_{\textrm{MA},a}(r,r)+\gamma_{\textrm{MA}}(r,r)\vert^2)\nonumber\\
&\leq  C \vert t-r\vert^{2\beta}.
\label{majoration d gamma}
\end{align}


The end of the proof is  a direct  application of Theorem 12.3 of   \cite{Billingsley1968}. 
Since $\beta>1/2,$ the sequence $n(\widehat{\gamma}_{\textrm{MA},a}(t,t)-\gamma_{\textrm{MA}}(t,t))$ is tight in $C([0,T])$ and converges in distribution to 0. The proof is complete with a direct application of the definition of weak convergence in $C([0,T])$ considering the bounded and continuous "sup" functional. 
 \hfill $\Box$
\subsection{Proofs related to the asymptotic normality and the confidence bands}
The steps of the proof of Proposition \ref{prop:CLTdansC} are similar to the steps of the proof of Proposition \ref{prop:gammaconv}. We first examine the finite combinations and invoke the Cramer-Wold device. Then we prove the tightness thanks to inequalities on the increments.

Let us first deal with the interpolation error, which is negligible under the assumption on the grid of discretization points, as shown in  (\ref{convergencemuinterpolee}).

Then, in light of (\ref{def:mutilde-muhat}), Lemma \ref{esperancecarrex} and Lemma \ref{espbetatinstantt}, we clearly have that, for each value of $t,$
\begin{align*}
\sqrt{n} \left( \widehat{\mu}_{\textrm{MA},a}(t)- \widetilde{\mu}(t) \right) &= o_p(1),
\end{align*}
and consequently, as $n$ tends to infinity,
\begin{align*}
\sqrt{n} \left( \widehat{\mu}_{\textrm{MA},a}(t)- \mu(t) \right) & \rightarrow  {\cal N}(0, \gamma_Z(t,t)) \quad \mbox{in distribution,}
\end{align*}
where the covariance-function of $\widetilde{\mu},$ which defined in (\ref{var_ma_asymp}), satisfies $\lim_{N \rightarrow \infty} n \gamma_{\textrm{MA}} = \gamma_Z.$

If we now consider $p$ distinct discretization instants $0\leq t_1 < t_2 \ldots < t_p \leq 1,$ it is immediate to check that for any vector ${\bf c} \in \mR^p,$ $\sqrt{n} \left(\sum_{j=1}^p c_j (\widetilde{\mu}(t_j) - \mu(t_j)) \right) \rightarrow {\cal N}(0,\sigma^2_c)$ where
\begin{align*}
\sigma^2_c &= \sum_{j=1}^p \sum_{\ell =1}^p c_j c_{\ell} \gamma_Z(t_j,t_{\ell}).
\end{align*}
Indeed, by linearity, there exists  a vector of random weights $(w_1, \ldots, w_N)$ which does not depend on time $t$ such that
\begin{align*}
\widetilde{\mu}(t) & =  \sum_{k \in U} w_k Y_k(t),
\end{align*}
and $\sum_{j=1}^p c_j \widetilde{\mu}(t_j) =  \sum_{k \in U} w_k \left( \sum_{j=1}^p c_j Y_k(t_j) \right)$ also satisfies a CLT, with asymptotic variance  $\sigma^2_c,$  under the moment conditions (A7). 
Thus, any finite linear combination is asymptotically Gaussian and we can conclude that the vector $\sqrt{n} \left(\widetilde{\mu}(t_1) - \mu(t_1), \ldots,  \widetilde{\mu}(t_p) - \mu(t_p) \right)$ is asymptotically Gaussian with the Cramer-Wold device.

It remains to check the tightness of the functional process and this  is a direct consequence of (\ref{distance mu hat - mu}) and (\ref{majoration d2}). Indeed, denoting by $Z_n(t) = \sqrt{n} \left( \widehat{\mu}_{\textrm{MA},a}(t)- \mu(t) \right),$ there is a constant $C$ such that,  for all $(r,t) \in [0,T]^2,$ 
\begin{align*}
\E_p \left( \left[ Z_n(t) - Z_n(r)\right]^2 \right) &\leq C \left| t-r \right|^{2\beta},
\end{align*}
and, since $\beta>1/2,$ the sequence $Z_n$ is tight in $C[0,T],$ in view of Theorem 12.3 of   \cite{Billingsley1968}. 

\hfill $\Box$

\medskip

We prove now Proposition \ref{conditional weak functional convergence for Gaussian processes}, the last result of the paper.
The proof consists in showing the weak convergence of the sequence of distributions $(\widehat{Z}_N)$ to the law of $Z$ in $C([0,T]).$ 

For any vector of $p$ points $0\leq t_1 < \ldots < t_p \leq T,$ the finite dimensional convergence of the distribution of the Gaussian vector $(\widehat{Z}_N(t_1),\ldots, \widehat{Z}_N(t_p))$ to the distribution of $(Z(t_1),\ldots, Z(t_p))$ is an immediate consequence of the uniform convergence of the covariance function stated in Proposition \ref{prop:gammaconv}. We can conclude with Slutsky's Lemma noting that for any $(c_1,\ldots, c_p) \in \mR^p,$ 
\begin{align}
 \sum_{j=1}^p \sum_{\ell =1}^p c_j c_{\ell} \widehat{\gamma}_{\textrm{MA},d}(t_j,t_{\ell}) & \rightarrow  \sum_{j=1}^p \sum_{\ell =1}^p c_j c_{\ell} {\gamma}_{\textrm{MA}}(t_j,t_{\ell}) \quad \mbox{in probability}.
\end{align}

Now, we need to check the tightness of $(\widehat{Z}_N)$ in $C([0,T]).$ Given $\widehat{\gamma}_{\textrm{MA},d},$ we have for $(r,t) \in [0,T]^2,$ 
\begin{align*}
\E_p \left[ \left( \widehat{Z}_N(t) - \widehat{Z}_N(r) \right)^2 | \ \widehat{\gamma}_{\textrm{MA},d} \right] & =  n \left( \widehat{\gamma}_{\textrm{MA},d}(t,t) - 2 \widehat{\gamma}_{\textrm{MA},d}(r,t) +  \widehat{\gamma}_{\textrm{MA},d}(r,r) \right) 
\end{align*} 
and after some algebra, we obtain thanks to  Assumption (A2) that
 \begin{align}
\E_p \left[ \left( \widehat{Z}_N(t) - \widehat{Z}_N(r) \right)^2 |  \widehat{\gamma}_{\textrm{MA},d} \right]   & \leq   \frac{C}{N}  \sum_{k \in U}  \left[ \left( Y_{k,d}(t) - Y_{k,d}(r)  \right)^2 + \left( \widehat{Y}_{k,d}(t) - \widehat{Y}_{k,d}(r)  \right)^2 \right]. \nonumber \\
  \label{def:ineqbootpar}
\end{align}

Let us first study the term $\sum_{k \in U}  \left( Y_{k,d}(t) - Y_{k,d}(r)  \right)^2$ in the previous inequality and without loss of generality suppose that $t>r.$ To check the continuity of the trajctories, we only need to consider points $r$ and $t$ that are close to each other.
If $t$ and $r$ belong to the same interval, say $[t_{i},t_{i+1}],$ then it  is easy to check, with Assumption (A4)  that
\begin{align}
\frac{1}{N} \sum_{k \in U} \left(Y_{k,d}(t) - Y_{k,d}(r) \right)^2 &=  \frac{(t-r)^2}{(t_{i+1}-t_i)^2} \frac{1}{N} \sum_{k \in U} \left(Y_k(t_{i+1}) - Y_k(t_{i}) \right)^2 \nonumber \\
  & \leq   C (t-r)^{2 \beta}.
\label{decomp:ykd2}  
\end{align}
If we suppose now that $r \in [t_{i-1},t_i]$ and $t \in [t_i, t_{i+1}],$ then we have
\begin{align*}
\frac{| Y_{k,d}(t) - Y_{k,d}(r)|}{t-r} &\leq  \max \left( \frac{|Y_k(t_{i+1}) - Y_k(t_{i})|}{t_{i+1}-t_{i}} , \frac{|Y_k(t_{i}) - Y_k(t_{i-1})|}{t_{i}-t_{i-1}}\right) \nonumber \\
 & \leq  \frac{|Y_k(t_{i+1}) - Y_k(t_{i})|}{t_{i+1}-t_{i}} + \frac{|Y_k(t_{i}) - Y_k(t_{i-1})|}{t_{i}-t_{i-1}}
\end{align*}
and using the same decomposition as in (\ref{decomp:ykd2}), we directly get that 
\begin{align*}
\sum_{k \in U}  \left( Y_{k,d}(t) - Y_{k,d}(r)  \right)^2 &\leq C (t-r)^{2 \beta}.
\end{align*}
The second term at the right-hand side of inequality (\ref{def:ineqbootpar}) is dealt with similar arguments and the decomposition used   in the proof of Lemma \ref{majoration norme beta}, so that $\frac{1}{N}  \sum_{k \in U}  \left( \widehat{Y}_{k,d}(t) - \widehat{Y}_{k,d}(r)  \right)^2 \leq C |t-r|^{2 \beta}.$  

Thus, the trajectories of the Gaussian process are continuous on $[0,T]$ whenever $\beta > 0$  (see {\it e.g} Theorem 1.4.1 in \cite{AdlerTaylor2007}) and the sequence $(\widehat{Z}_N)$ converges weakly to $Z$ in $C([0,T])$ equipped with the supremum norm.

Using again Proposition \ref{prop:gammaconv}, we have, uniformly in $t,$  $\widehat{\sigma}_Z(t) = \sigma_Z(t) + o_p(1),$ where  $\widehat{\sigma}_Z^2(t)  = n \widehat{\gamma}_{\textrm{MA},d}(t,t).$ Since, by hypothesis $\sigma^2_Z(t) =\gamma_Z(t,t)$ is a continuous function and $\inf_t \gamma_Z(t,t) >0,$ we get with Slutsky's lemma that $(\widehat{Z}_N/ \widehat{\sigma}_Z)$ converges weakly to $Z/\sigma_Z$ in $C([0,T]).$ By definition of the weak convergence in $C([0,T])$ and the continuous mapping theorem, we also deduce that the real random variable $\widehat{M}_N = \sup_{t \in [0,T]} |\widehat{Z}_N(t)|/ \widehat{\sigma}_Z(t)$ converges in distribution to $M = \sup_{t \in [0,T]} |Z(t)|/ \sigma_Z(t),$ so that for each $c\geq 0$, 
\begin{align*}
\P \left(\sup_{t \in [0,T]} |\widehat{Z}_N(t)|/\widehat{\sigma}_Z(t) \leq c \right) & \rightarrow  \P \left(\sup_{t \in [0,T]} | Z(t)| /\sigma_Z(t) \leq c \right).
\end{align*}

Note finally, that under the previous  hypotheses on $\gamma_Z$ (see {\it e.g.} \cite{PittTran1979}), the real random variable $M= \sup_{t \in [0,T]} \left(  | Z(t)| / \sigma_Z(t) \right)$ has an absolutely continuous and bounded density function so that the convergence holds uniformly in $c$ (see {\it e.g.} Lemma  2.11 in \cite{vanderVaart1998}).
\hfill $\Box$

\subsection{Some  useful lemmas}
We state here without any proof some results  that are needed for the study of the convergence of the covariance function. They rely on applications of the Cauchy-Schwarz inequality and on  the assumptions on the moments of the trajectories and the inclusion probabilities.
\begin{Lem}
Assume (A2)-(A5) and (A7) hold.  There are two constants  $\zeta_4$ and $\zeta_5$ such that
$$
\frac{1}{N}\sum_{k\in U}\tilde{e}_k(t)^2\tilde{e}_k(r)^2\leq\zeta_4
$$
and 
$$
\frac{1}{N^2}\sum_{k\in U}\sum_{l\in U}\tilde{e}_k(t)^2\tilde{e}_l(r)^2\leq\zeta_5,
$$
where  $\tilde{e}_k(t)=Y_k(t)-\tilde{Y}_k(t).$
\label{somme e tilde puissance 2}
\end{Lem}

\begin{Lem}
Assume (A2)-(A5) and (A7) hold. There are  two constants $\zeta_{6}$ and $\zeta_{7}$ such that
\begin{equation*}
\E_p\left(\frac{1}{N}\sum_{k\in U}\widehat{\tilde{\phi}}_{kk}(t,r)^2\right)\leq \zeta_{6}\vert t-r\vert^{2\beta}
\end{equation*}
and
\begin{equation*}
\E_p\left(\frac{1}{N^2}\sum_{k,l\in U}\widehat{\tilde{\phi}}_{kl}(t,r)\right)^2\leq \zeta_{7}[\vert t-r\vert^{2\beta}
\end{equation*}
where $\widehat{\tilde{\phi}}_{kl}(t,r)=\widehat{\tilde{e}}_k(t)\widehat{\tilde{e}}_l(t)-\widehat{\tilde{e}}_k(r)\widehat{\tilde{e}}_l(r)$ and $\widehat{\tilde{e}}_k(t)=\tilde{Y}_k(t)-\widehat{Y}_{k,a}(t)$.
\label{somme phi tilde hat puissance 2}
\end{Lem}

\begin{Lem}
Assume (A2)-(A5) and (A7) hold. There are two constant constants  $\zeta_8$ and $\zeta_{9}$ such that 
$$
\frac{1}{N}\sum_{k\in U}\phi_{kk}^2(t,r)\leq \zeta_8\vert t-r\vert^{2\beta}
$$
and 
$$
\left(\frac{1}{N^2}\sum_{k,l\in U}\phi_{kl}(t,r)\right)^2\leq \zeta_{9}\vert t-r\vert^{2\beta}
$$
\label{somme phi puissance 2}
where $\phi_{kl}(t,r)=\tilde{e}_k(t)\tilde{e}_l(t)-\tilde{e}_k(r)\tilde{e}_l(r)$ and $\tilde{e}_k(t)=Y_k(t)-\tilde{Y}_k(t)$.
\end{Lem}

\begin{Lem}
Assume (A2)-(A5) and (A7) hold. There are two constants $\zeta_{10}$ and $\zeta_{11}$such that
\begin{equation*}
\E_p\left(\frac{1}{N}\sum_{k\in U}\tilde{\phi}_{kk}(t,r)^2\right)\leq \zeta_{10}\vert t-r\vert^{2\beta}
\end{equation*}
and
\begin{equation*}
\E_p\left(\frac{1}{N^2}\sum_{k,l\in U}\tilde{\phi}_{kl}(t,r)\right)^2\leq \zeta_{11}\vert t-r\vert^{2\beta}
\end{equation*}
\label{somme phi tilde puissance 2}
where $\tilde{\phi}_{kl}(t,r)=\tilde{e}_k(t)\widehat{\tilde{e}}_l(t)-\tilde{e}_k(r)\widehat{\tilde{e}}_l(r)$, $\tilde{e}_k(t)=Y_k(t)-\tilde{Y}_k(t)$ and $\widehat{\tilde{e}}_k(t)=\tilde{Y}_k(t)-\widehat{Y}_{k,a}(t)$.
\end{Lem}

\bibliographystyle{apalike}


\begin{thebibliography}{}

\bibitem[Adler and Taylor, 2007]{AdlerTaylor2007}
Adler, R.~J. and Taylor, J.~E. (2007).
\newblock {\em Random Fields and Geometry}.
\newblock Springer-Verlag, New York.

\bibitem[Bhatia, 1997]{Bhatia96}
Bhatia, R. (1997).
\newblock {\em Matrix Analysis}.
\newblock Springer-Verlag, New York.

\bibitem[Billingsley, 1968]{Billingsley1968}
Billingsley, P. (1968).
\newblock {\em Convergence of Probability Measures}.
\newblock John Wiley and Sons.

\bibitem[Boistard et~al., 2012]{BLRG2012}
Boistard, H., Lopuha{\"a}, H.~P., and Ruiz-Gazen, A. (2012).
\newblock Approximation of rejective sampling inclusion probabilites and
  application to higher order correlation.
\newblock {\em Electronic J. of Statistics}, 6:1967--1983.


\bibitem[Bosq, 2000]{Bosq_2000}
Bosq, D. (2000).
\newblock {\em Linear Processes in Function Spaces: Theory and Applications},
  volume 149 of {\em Lecture notes in Statistics}.
\newblock Springer-Verlag, New York.

\bibitem[Breidt and Opsomer, 2000]{breidt_opsomer_2000}
Breidt, F. J. and Opsomer, J. D.(2000).
\newblock Local polynomial regression estimators in survey sampling.
\newblock {\em Ann.  Statist.}, 28(4):1023--1053.

\bibitem[Callado et~al., 2009]{callado2009}
Callado, A., Kamienski, C., Szab\'o, G.,  Ger\"o, B., Kelner, J., Fernandes, S. and  Sadok, D. (2009).
\newblock A Survey on Internet Traffic Identification and Classification.
\newblock {\em IEEE Communications Surveys and Tutorials}, 11:37--52.

\bibitem[Cardot, 2007]{Cardot2007}
Cardot, H. (2007).
\newblock Conditional functional principal components analysis.
\newblock {\em Scandinavian J. of Statistics}, 34:317--335.

\bibitem[Cardot et~al., 2010]{cardot_chaouch_goga_labruere_2010}
Cardot, H., Chaouch, M., Goga, C., and Labru\`ere, C. (2010).
\newblock Properties of design-based functional principal components analysis.
\newblock {\em J. of Statistical Planning and Inference}, 140:75--91.

\bibitem[Cardot et~al., 2012a]{CardotDegrasJosserand2012}
Cardot, H., Degras, D., and Josserand, E. (2012a).
\newblock Confidence bands for {H}orvitz-{T}hompson estimators using sampled
  noisy functional data.
\newblock  To appear in {\em Bernoulli}.

\bibitem[Cardot et~al., 2012b]{CDGJL2012}
Cardot, H., Dessertaine, A., Goga, C., Josserand, E., and Lardin, P. (2012b).
\newblock Comparaison de diff\'erents plans de sondage et construction de
  bandes de confiance pour l'estimation de la moyenne de donn\'ees
  fonctionnelles : une illustration sur la consommation \'electrique.
\newblock To appear in {\em Techniques d'Enqu\^etes/Survey Methodology}.

\bibitem[Cardot et~al., 2012c]{CGL2012}
Cardot, H., Goga, C.,  and Lardin, P. (2012c).
\newblock Variance estimation and asymptotic confidence bands for the mean estimator of sampled functional data with high entropy unequal probability sampling designs.
\newblock {\em Arxiv:1209.6503}

\bibitem[Cardot and Josserand, 2011]{cardot_josserand_2011}
Cardot, H. and Josserand, E. (2011).
\newblock Horvitz-{T}hompson estimators for functional data: asymptotic
  confidence bands and optimal allocation for stratified sampling.
\newblock {\em Biometrika}, 98:107--118.

\bibitem[Chiou et~al., 2004]{CMW2004}
Chiou, J., M\"uller, H., and Wang, J. (2004).
\newblock Functional response models.
\newblock {\em Statistica Sinica}, 14:675--693.

\bibitem[Cuevas et~al., 2006]{CuevasFF06}
Cuevas, A., Febrero, M., and Fraiman, R. (2006).
\newblock On the use of the bootstrap for estimating functions with functional
  data.
\newblock {\em Computational Statistics and Data Analysis}, 51:1063--1074.

\bibitem[Degras, 2011]{degras_2011}
Degras, D. (2011).
\newblock Simultaneous confidence bands for parametric regression with
  functional data.
\newblock {\em Statistica Sinica}, 21(4):1735--1765.

\bibitem[Degras, 2012]{Degras2012}
Degras, D. (2012).
\newblock Rotation sampling for functional data.
\newblock http://arxiv.org/abs/1204.4494.

\bibitem[Deville and S\"arndal, 1992]{DevilleSarndal1992}
Deville, J.~C. and S\"arndal, C.~E. (1992).
\newblock Calibration estimators in survey sampling.
\newblock {\em J. Amer. Statist. Assoc.}, 87:376--382.

\bibitem[Faraway, 1997]{faraway_1997}
Faraway, J. (1997).
\newblock Regression analysis for a functional response.
\newblock {\em Technometrics}, 39(3):254--261.

\bibitem[Ferraty et~al., 2011]{FLTV2011}
Ferraty, F., Laksaci, A., Tadj, A., and Vieu, P. (2011).
\newblock Kernel regression with functional response.
\newblock {\em Electronic J. of Statist.}, 5:159--171.

\bibitem[Fuller, 2009]{FullerLivre}
Fuller, W.~A. (2009).
\newblock {\em Sampling Statistics}.
\newblock John Wiley and Sons, Hoboken, New Jersey.

\bibitem[Guillas, 2001]{Guillas2001}
Guillas, S. (2001).
\newblock Rates of convergence of autocorrelation estimates for autoregressive
  {H}ilbertian processes.
\newblock {\em Statist. and Probability Letters}, 55:281--291.

\bibitem[Hahn, 1977]{Hahn77}
Hahn, M. (1977).
\newblock Conditions for sample-continuity and the central limit theorem.
\newblock {\em Annals of Probability}, 5:351--360.

\bibitem[H\'ajek, 1981]{Hajek1981}
H\'ajek, J. (1981).
\newblock {\em Sampling From a Finite Population}.
\newblock Statistics: Textbooks and Monographs. Marcel Dekker, New York.

\bibitem[Isaki and Fuller, 1982]{isaky_fuller_1982}
Isaki, C. and Fuller, W. (1982).
\newblock Survey design under the regression superpopulation model.
\newblock {\em J. Amer. Statist. Assoc.}, 77:49--61.

\bibitem[Pitt and Tran, 1979]{PittTran1979}
Pitt, L.~D. and Tran, L.~T. (1979).
\newblock Local sample path properties of {G}aussian fields.
\newblock {\em Annals of Probability}, 7:477--493.

\bibitem[Ramsay and Silverman, 2005]{Ramsay_Silverman_Livre}
Ramsay, J.~O. and Silverman, B.~W. (2005).
\newblock {\em Functional Data Analysis}.
\newblock Springer-Verlag, New York, second edition.

\bibitem[Robinson and S\"arndal, 1983]{RobinsonSarndal83}
Robinson, P. and S\"arndal, C. (1983).
\newblock Asymptotic properties of the generalized regression estimator in
  probability sampling.
\newblock {\em Sankhya : The Indian Journal of Statistics}, 45:240--248.

\bibitem[S\"arndal, 1980]{Sarndal1980}
S\"arndal, C. (1980).
\newblock On $\pi$ inverse weighting versus best linear unbiased weighting in
  probability sampling.
\newblock {\em Biometrika}, 67:639--50.

\bibitem[S\"{a}rndal et~al., 1992]{SarndalLivre}
S\"{a}rndal, C.~E., Swensson, B., and Wretman, J. (1992).
\newblock {\em Model Assisted Survey Sampling}.
\newblock Springer series in statistics. Springer-Verlag, New York.

\bibitem[van~der Vaart, 1998]{vanderVaart1998}
van~der Vaart, A.~W. (1998).
\newblock {\em Asymptotic Statistics}.
\newblock Cambridge University Press.

\bibitem[van~der Vaart and Wellner, 2000]{vdVWellener2000}
van~der Vaart, A.~W. and Wellner, J.~A. (2000).
\newblock {\em Weak Convergence and Empirical Processes. With Applications to
  Statistics}.
\newblock Springer-Verlag, New York.

\end{thebibliography}

\end{document}